\documentclass[10pt]{amsart}

\usepackage{amsmath}
\usepackage{amsthm}
\usepackage{hyperref}

\newtheorem{prop}{Proposition}[section]
\newtheorem{lemma}[prop]{Lemma}
\newtheorem{thm}[prop]{Theorem}
\newtheorem{cor}[prop]{Corollary}

\theoremstyle{definition}
\newtheorem{defn}[prop]{Definition}

\newcommand{\nm}[2]{|| #1 ||_{#2}}
\newcommand{\dt}{\frac{\partial}{\partial t}}
\newcommand{\brs}[1]{\left| #1 \right|}

\newcommand{\gG}{\Gamma}

\newcommand{\gD}{\Delta}
\newcommand{\gd}{\delta}
\newcommand{\gl}{\lambda}
\newcommand{\gt}{\theta}
\newcommand{\gw}{\omega}
\newcommand{\ga}{\alpha}
\newcommand{\gb}{\beta}
\renewcommand{\ge}{\epsilon}
\newcommand{\N}{\nabla}
\newcommand{\FF}{\mathcal F}
\renewcommand{\bar}[1]{\overline{#1}}
\newcommand{\del}{\partial}
\newcommand{\delb}{\bar{\partial}}
\newcommand{\bi}{\bar{i}}
\newcommand{\bj}{\bar{j}}
\newcommand{\bk}{\bar{k}}
\newcommand{\bl}{\bar{l}}
\newcommand{\bm}{\bar{m}}
\newcommand{\bp}{\bar{p}}
\newcommand{\bn}{\bar{n}}
\newcommand{\bq}{\bar{q}}
\newcommand{\br}{\bar{r}}
\newcommand{\bs}{\bar{s}}

\newcommand{\bv}{\bar{v}}
\newcommand{\til}[1]{\widetilde{#1}}
\newcommand{\sq}{\square}
\newcommand{\ohat}[1]{\overset{\circ}{#1}}

\newcommand{\hook}{\mathbin{\hbox{\vrule height2.4pt width4.5pt depth-2pt
\vrule height5pt width0.4pt depth-2pt}}}

\newcommand{\static}{static}

\DeclareMathOperator{\Sym}{Sym}
\DeclareMathOperator{\divg}{div}
\DeclareMathOperator{\Rc}{Rc}
\DeclareMathOperator{\Rm}{Rm}
\DeclareMathOperator{\tr}{tr}

\begin{document}

\title{Hermitian Curvature Flow}

\author{Jeffrey Streets}
\address{Fine Hall\\
         Princeton University\\
         Princeton, NJ 08544}
\email{\href{mailto:jstreets@math.princeton.edu}{jstreets@math.princeton.edu}}

\author{Gang Tian}
\address{Fine Hall\\
	 Princeton University\\
	 Princeton, NJ 08544}
\email{\href{mailto:tian@math.princeton.edu}{tian@math.princeton.edu}}

\thanks{The first author was supported by the National Science Foundation via
DMS-0703660}
\thanks{The second author was partly supported by the National Science
Foundation via
DMS-0703985 and DMS-0804095}

\date{January 26, 2009}

\begin{abstract} We define a functional for Hermitian metrics using the
curvature of the Chern connection.  The Euler-Lagrange equation for this
functional is an elliptic equation for Hermitian metrics.  Solutions to this
equation are related to K\"ahler-Einstein metrics, and are automatically
K\"ahler-Einstein under certain conditions.  Given this, a natural
parabolic flow equation arises.  We prove short time existence and regularity
results for this flow, as well as stability for the flow near K\"ahler-Einstein
metrics with negative or zero first Chern class.
\end{abstract}

\maketitle

\section{Introduction}

In this paper we introduce a new curvature evolution equation on compact complex
manifolds.
Specifically, given $(M^{2n}, g, J)$ a manifold with integrable
complex structure $J$ and Hermitian metric $g$, let $\N$ denote the Chern
connection of $g$, which is a metric compatible connection with torsion $T$
\cite{Morrow}.  Let
$\Omega$ denote the curvature of $\N$.  Define
\begin{gather*}
S_{i \bj} = \left(\tr_{\gw} \Omega \right)_{i \bj} = g^{k \bl}
\Omega_{k \bl i \bj}
\end{gather*}
and let $s = g^{i \bj} S_{i \bj}$ be the scalar Chern curvature.  Furthermore
let $w_i = g^{j \bk} T_{i j \bk}$ denote the trace of the torsion.  Consider the
functional
\begin{align} \label{functional}
\mathbb F(g) =&\ \frac{\int_M \left[ s - \frac{1}{4} \brs{T}^2 - \frac{1}{2}
\brs{w}^2 \right] dV}{\left(\int_M dV \right)^{\frac{n-1}{n}}}.
\end{align}
As we will see in section 3 this is the unique functional yielding $\ohat{S}$ as
the traceless component of the second-order terms in the associated
Euler-Lagrange equation.  Moreover, the form of the Euler-Lagrange equation
suggests a flow equation in the same way
that Ricci flow is suggested by the usual Hilbert functional.  In particular we
define an evolution equation
\begin{gather} \label{HCF}
\frac{\del}{\del t} g = - S + Q
\end{gather}
where $Q = Q(T)$ is a certain quadratic polynomial in the torsion $T$ of $\N$
which
is made precise in section 3.  We call equation (\ref{HCF}) \emph{Hermitian
curvature flow} (HCF).  Of course now it is known that Ricci flow
is indeed the gradient flow of the lowest eigenvalue of a certain Schr\"odinger
operator, although a corresponding statement for HCF is not yet known.  It is
also possible to
write HCF in terms of Hodge-type operators.  In particular, if $\gw(t)$ denotes
the K\"ahler form of the time varying metric, then it satisfies the equation
\begin{align*}
\frac{\del}{\del t} \gw =&\ - \left( \del_g^* \del \gw - \del \del_g^* \gw -
\frac{\sqrt{-1}}{2} \del \delb \log \det g - 2 \sqrt{-1} \left( \delb^*_g \gw
\hook \delb \gw \right) \right) + Q'.
\end{align*}
where here $Q'$ is a distinct fixed quadratic expression in the torsion.

We observe that when a solution $g(t)$ to HCF exists, the metric
is Hermitian with respect
to the fixed complex structure $J$ for all time.  Secondly, we will show that
when the initial metric $g(0)$
is K\"ahler, then the solution $g(t)$ is K\"ahlerian and consequently the
solution to HCF is
given by K\"ahler-Ricci flow.
Thirdly, we prove that certain static solutions are K\"ahler-Einstein metrics.
It will be a very interesting problem to classify all static solutions. It is
possible that they are all K\"ahler.
Hence, in some sense, this new flow evolves Hermitian metrics towards
K\"ahler metrics.

Next we will prove a local existence
theorem for HCF and develop some regularity properties for this flow.  In
particular we derive higher order derivative estimates in the presence of a
curvature bound.  A consequence of
these estimates is the following short-time existence theorem.
\begin{thm} \label{STEThm} Let $(M^{2n}, g_0, J)$ be a complex manifold with
Hermitian metric $g_0$.  There exists a constant $c(n)$ depending only on the
dimension such that there exists a unique solution $g(t)$ to HCF for
\begin{align*}
t \in \left[0, \frac{c(n)}{\max \{ \brs{\Omega}_{C^0(g_0)}, \brs{\N
T}_{C^0(g_0)}, \brs{T}^2_{C^0(g_0)}
\}} \right].
\end{align*}
Moreover, there exist constants $C_m$ depending only on $m$ such that the
estimates
\begin{align*}
\brs{\N^m \Omega}_{C^0(g_t)}, \brs{\N^{m+1} T}_{C^0(g_t)} \leq \frac{C_m \max \{
\brs{\Omega}_{C^0(g_0)}, \brs{\N T}_{C^0(g_0)},
\brs{T}^2_{C^0(g_0)} \}}{t^{m / 2}}
\end{align*}
hold for all $t$ in the above interval.  Moreover, the solution exists on a
maximal time interval $[0, \tau)$, and if $\tau < \infty$ then
\begin{align*}
\limsup_{t \to \tau} \max \{ \brs{\Omega}_{C^0(g_t)}, \brs{\N T}_{C^0(g_t)},
\brs{T}_{C^0(g_t)} \} = \infty.
\end{align*}
\end{thm}

In some sense, the simplest possible behavior for this flow should occur near
K\"ahler-Einstein metrics, where we expect the flow to be not too much different
from K\"ahler-Ricci flow.  To that end we prove a stability result for HCF
around K\"ahler-Einstein metrics with negative or zero first Chern class.
Specifically, we show
\begin{thm} \label{Stabthm} Let $(M^{2n}, g, J)$ be a complex manifold with
K\"ahler-Einstein
metric $g$ and $c_1(M) < 0$ or $c_1(M) = 0$.  Then there exists $\ge = \ge(g)$
so that if
$\til{g}$ is a Hermitian metric on $M$
compatible with $J$ and $\brs{\til{g} - g}_{C^{\infty}} < \ge$ then the solution
to HCF with initial condition $\til{g}$ exists for all time and converges to a
K\"ahler-Einstein metric.
\end{thm}

There are two natural directions which motivate defining this flow.  First,
given all of the success of Ricci flow it is natural to study it on complex
manifolds.  However, it is usually the case that the Ricci tensor of a Hermitian
metric is not $(1,1)$, and thus the Hermitian condition for the metric is not
preserved.  Thus
the Ricci flow is not the best tool for studying complex geometry which is not
already K\"ahler.  The tensor $S$ is a natural $(1,1)$ curvature tensor
associated to a Hermitian metric which differs from the Ricci tensor by torsion
terms, meaning that it equals the Ricci tensor in the K\"ahler setting.
Moreover, the operator $g \to S(g)$ is strictly elliptic, giving HCF nice
existence properties.  Thus from this perspective HCF is the right analogue of
Ricci flow for Hermitian geometry.

The second motivation, and actually our original motivation for HCF is that it
serves as a ``holonomy flow.''  If one looks on the level of the K\"ahler form
and asks for a parabolic flow which preserves the Hermitian condition and is
stationary on K\"ahler manifolds, HCF comes up quite naturally.  There is the
side effect that one ends up looking not just for K\"ahler metrics,
but \emph{K\"ahler-Einstein} metrics.  Given the excellent existence properties
of the K\"ahler-Ricci flow, this is an acceptable price to pay.  Indeed, other
natural analytic approaches to this question which strictly look for K\"ahler
metrics among Hermitian metrics (see for instance \cite{Vaisman}) yield
equations which are not elliptic.  In \cite{Catenacci} it is shown that if the
usual Ricci-type curvature of the Chern connection is a nonzero scalar multiple
of the metric, then the metric is automatically K\"ahler-Einstein.  However,
this Ricci tensor is not in general $(1,1)$, so from the perspective of
Hermitian geometry, especially defining a flow of Hermitian metrics, this
condition is not natural.

Here is an outline of the rest of the paper.  In section 2 we define
all of
the relevant objects and notation and provide various curvature formulas.  In
section 3 we discuss the Hermitian Hilbert functional.  Section 4 gives the
definition of HCF and provides various equivalent formulations using Hodge-type
operators and the Levi-Civita connection.  In sections 5-7 we prove existence
and regularity properties for HCF.  In section 8 we prove the stability result
for HCF around K\"ahler-Einstein metrics.  We conclude in section 9 with a
discussion of some related questions.  Section 10 is an appendix containing
various useful calculations related to Hermitian geometry.

\textbf{Acknowledgments} The first author would like to thank Aaron Naber,
Yanir Rubinstein, and Jeff Viaclovsky for interesting conversations.

\section{Differential Operators on Hermitian Manifolds}
Let $(M^{2n}, g, J)$ be a complex manifold with a Hermitian metric $g$.  In
particular
\begin{gather*}
J: TM \to TM
\end{gather*}
is an integrable almost complex structure, i.e.
\begin{gather*}
N_J(X,Y) := [JX, JY] - J[JX, Y] - J[X,JY] - [X,Y] = 0
\end{gather*}
for all $X, Y \in TM_p$.  Furthermore
\begin{gather*}
g(u,v) = g(Ju, Jv).
\end{gather*}
This equation is written in a unitary frame as
\begin{gather*}
g_{i j} = g_{\bi \bj} = 0, \qquad g_{i \bj} = g_{\bj i} = \bar{g}_{\bi j}.
\end{gather*}
First we recall the Chern connection $\N$.  In complex coordinates, the only
nonvanishing components of the connection are given by
\begin{align*}
\gG_{i j}^k =&\ g^{k \bl} \del_i g_{j \bk}.
\end{align*}
This connection is compatible with $g$, but has torsion $T$.  In particular in
complex coordinates we have
\begin{align*}
T_{i j}^k =&\ g^{k \bl} \left( \del_i g_{j \bl} - \del_j g_{i \bl} \right).
\end{align*}
Also, there is a natural trace of the torsion
\begin{gather} \label{wdef}
w_i = T_{i j}^j.
\end{gather}
As we will see in the next section $w$ is just a multiple of $\delb^* \gw$ but
this separate definition will be useful to us.  We will also need certain
quadratic expressions in the torsion.  Specifically let
\begin{gather} \label{Qdefs}
\begin{split}
Q^1_{i \bj} =&\ g^{k \bl} g^{m \bn} T_{i k \bn} T_{\bj \bl m}\\
Q^2_{i \bj} =&\ g^{k \bl} g^{m \bn} T_{\bl \bn i} T_{k m \bj}\\
Q^3_{i \bj} =&\ g^{k \bl} g^{m \bn} T_{i k \bl} T_{\bj \bn m}\\
Q^4_{i \bj} =&\ \frac{1}{2} g^{k \bl} g^{m \bn} \left(T_{m k \bl} T_{\bn \bj i}
+ T_{\bn \bl
k} T_{m i \bj} \right).
\end{split}
\end{gather}
Note that each $Q^i$ is a real symmetric $(1,1)$ tensor.  The covariant
derivatives of torsion also satisfy an identity.
\begin{lemma} \label{Tformula} Given $g$ a Hermitian metric,
\begin{align*}
\N_i T_{j k \bl} + \N_k T_{i j \bl} + \N_j T_{k i \bl} =&\ T_{i j}^p T_{k p \bl}
+ T_{j k}^p T_{i p \bl} + T_{k i}^p T_{j p \bl}
\end{align*}
\begin{proof} We directly compute
\begin{align*}
\N_i T_{j k \bl} + \N_k T_{i j \bl} + \N_j T_{k i \bl} =&\ \del_i T_{j k \bl} +
\del_k T_{i j \bl} + \del_j T_{k i \bl}\\
&\ - \gG_{i j}^p T_{p k \bl} - \gG_{i k}^p T_{j p \bl} - \gG_{k i}^p T_{p j
\bl}\\
&\ - \gG_{k j}^p T_{i p \bl} - \gG_{j k}^p T_{p i \bl} - \gG_{j i}^p T_{k p
\bl}\\
=&\ T_{i j}^p T_{k p \bl} + T_{j k}^p T_{i p \bl} + T_{k i}^p T_{j p \bl}
\end{align*}
\end{proof}
\end{lemma}
Next we collect some useful formulas for the Chern curvature.  In particular,
let $\Omega$
denote the curvature of the Chern connection and let $S$ be the trace, i.e.
\begin{gather*}
S_{\ga \bar{\gb}} = \left(\tr_{\gw} \Omega \right)_{\ga \bar{\gb}} = g^{\mu
\bar{\nu}} \Omega_{\mu \bar{\nu} \ga \bar{\gb}}.
\end{gather*}
Further, let
\begin{gather*}
s = g^{\ga \bar{\gb}} S_{\ga \bar{\gb}}.
\end{gather*}
We will let $P$ denote the trace of the transpose of $\Omega$, namely
\begin{align*}
P_{\ga \bar{\gb}} = g^{\mu \bar{\nu}} \Omega_{\ga \bar{\gb} \mu \bar{\nu}}.
\end{align*}
In the K\"ahler case $S_{\ga \bar{\gb}} = P_{\ga \bar{\gb}}$ is the Ricci
curvature and $s = r$ is the
scalar curvature.
\begin{lemma} \label{Scalc} Given $g$ a Hermitian metric we have
\begin{gather*}
\begin{split}
S_{j \bk} =&\ - g^{l \bm} g_{j \bk, l \bm} + g^{l \bm} g^{p \bq} g_{p \bk, \bm}
g_{j \bq,
l}
\end{split}
\end{gather*}
\begin{proof} First of all we have
\begin{align*}
\Omega_{l \bm j \bk} =&\ - g_{\bk p} \del_{\bm} \left( g^{p \bq} \del_{l} g_{j
\bq}
\right)\\
=&\ - g_{j \bk, l \bm} + g^{p \bq} g_{p \bk, \bm} g_{j \bq, l}
\end{align*}
Thus
\begin{align*}
S_{j \bk} =&\ g^{l \bm} \Omega_{l \bm j \bk}\\
=&\ - g^{l \bm} g_{j \bk, l \bm} + g^{l \bm} g^{p \bq} g_{p \bk, \bm} g_{j \bq,
l}
\end{align*}
as required.
\end{proof}
\end{lemma}

\begin{lemma} \label{Bianchi} \emph{(Bianchi Identity)} For $X, Y, Z \in T_x(M)$
we have
\begin{align*}
\Sigma \{ \Omega(X, Y) Z \} =&\ \Sigma \{ T(T(X, Y), Z) + \N_X T(Y, Z) \}\\
\Sigma \{ \N_X \Omega(Y, Z) + \Omega(T(X, Y), Z) \} =&\ 0
\end{align*}
\end{lemma}

\begin{lemma} \label{PSForm} Given $g$ a Hermitian metric we have
\begin{align*}
P_{i \bj} - S_{i \bj} =&\ g^{k \bl} \left( \N_{\bl} T_{k i \bj} + \N_i T_{\bl
\bj k} \right)
\end{align*}
\begin{proof} We compute using the Bianchi identity and the symmetries of the
torsion
\begin{align*}
\Omega_{i \bj k \bl} =&\ \Omega_{k \bj i \bl} + \N_{\bj} T_{k i \bl}\\
=&\ \Omega_{\bj k \bl i} + \N_{\bj} T_{k i \bl}\\
=&\ \Omega_{\bl k \bj i} + \N_k T_{\bl \bj i} + \N_{\bj} T_{k i \bl}.
\end{align*}
Taking the trace and relabelling indices gives the result.
\end{proof}
\end{lemma}

Now we focus on Hodge operators associated to $g$.  Let
\begin{gather*}
\gw(u,v) = - g(u, Jv)
\end{gather*}
be the K\"ahler form of $g$.  In local complex coordinates we have
\begin{gather*}
\gw = \frac{\sqrt{-1}}{2} g_{i \bj} dz^i \wedge d \bar{z}^j.
\end{gather*}
Let
\begin{gather*}
\Lambda^k = \bigoplus_{p+q = k} \Lambda^{p,q}
\end{gather*}
denote the usual decomposition of complex differential two-forms into forms of
type $(p,q)$.  The exterior differential $d$ decomposes into the operators
$\del$ and $\delb$
\begin{align*}
\del &: \Lambda^{p,q} \to \Lambda^{p+1, q}\\
\delb &: \Lambda^{p,q} \to \Lambda^{p, q+1}.
\end{align*}
Also the operator $d^*_{g}$, the $L^2$ adjoint of $d$, decomposes into
$\del^*_{g}$ and $\delb^*_{g}$
\begin{align*}
\del^*_g :&\ \Lambda^{p+1}, q \to \Lambda^{p, q}\\
\delb^*_g :&\ \Lambda^{p, q+1} \to \Lambda^{p,q}
\end{align*}
Using these operators we can define the complex Laplacians
\begin{align*}
\sq_\gw =&\ \del^*_g \del + \del \del^*_g : \Lambda^{p,q} \to
\Lambda^{p,q}\\
\bar{\sq}_\gw =&\ \delb^*_g \delb + \delb \delb^*_g : \Lambda^{p,q} \to
\Lambda^{p,q}
\end{align*}
It is well known that the operator $\ga \to \sq_\gw \ga$ is a second-order
elliptic operator with symbol that of the Laplacian in complex coordinates
\cite{Morrow}.  Moreover, one has the formula
\begin{gather} \label{LaplacianDecomposition}
\gD_{d, g} = \sq_{g} + \bar{\sq}_{g} + \mbox{lower order terms}
\end{gather}
However, we will be interested in the action of these operators on $\gw$ itself,
so the terms which are lower-order in (\ref{LaplacianDecomposition}) become
highest order terms in this context.  In the lemmas which follow we compute the
action of these differential operators explicitly.

\begin{lemma} \label{operatorlemma0} Given $g$ a Hermitian metric we have in
complex coordinates
\begin{align}
\left(\del^*_{g} \gw \right)_{\bk} =&\ \frac{\sqrt{-1}}{2} g^{p \bq}
\left(\del_{\bq} g_{p
\bk} - \del_{\bk} g_{p \bq} \right)\\
\left( \delb_g^* \gw \right)_j =&\ \frac{\sqrt{-1}}{2} g^{p \bq} \left( \del_p
g_{j \bq} - \del_j g_{p \bq} \right)
\end{align}
\begin{proof} We compute using integration by parts.  Given $\ga \in
\Lambda^{0,1}$ we have
\begin{align*}
\left( \del_{g}^* \gw, \ga \right) =&\ \left( \gw, \del \ga \right)\\
=&\ \int_M g^{\bk l} g^{\bi j} \left(\gw_{j \bk} \bar{ \del \ga_{i \bl}} \right)
\bar{g}\\
=&\ \frac{\sqrt{-1}}{2} \int_M g^{\bi l} \left( \bar{\ga_{\bl, i}} \right)
\bar{g}\\
=&\ - \frac{\sqrt{-1}}{2} \int_M \bar{\ga_{\bl}} \left[ \del_{\bi} \left( g^{\bi
l} \bar{g} \right) \right]\\
=&\ - \frac{\sqrt{-1}}{2} \int_M \bar{\ga_{\bl}} \left(\bar{g}\right) \left[ -
g^{\bi m} \del_{\bi}
g_{m \bn} g^{\bn l} + g^{\bi l} \frac{1}{\bar{g}} \del_{\bi} \bar{g} \right].
\end{align*}
This gives the first formula, and the second follows analogously.
\end{proof}
\end{lemma}

\begin{lemma} \label{operatorlemma1} Given $g$ a Hermitian metric we have in
complex coordinates
\begin{align}
\left( \del \del^*_{g} \gw \right)_{j \bk} =&\ \frac{\sqrt{-1}}{2} \left[g^{p
\bq} \left( g_{p \bk, \bq j}
- g_{p \bq, \bk j} \right) + g^{p \bq} g^{r \bs} g_{r \bq, j} \left(  g_{p \bs,
\bk} - g_{p \bk, \bs}\right) \right]
\end{align}
\begin{proof} In general for $\ga \in \Lambda^{0, 1}$ we have
\begin{align*}
\left(\del \ga \right)_{j \bk} = \del_j \ga_{\bk}.
\end{align*}
Thus we compute using Lemma \ref{operatorlemma0}
\begin{align*}
\left(\del \del_{g}^* \gw \right)_{j \bk} =&\ \frac{\sqrt{-1}}{2} \del_j \left(
g^{p \bq}
\left(\del_{\bq} g_{p \bk} - \del_{\bk} g_{p \bq} \right) \right)\\
=&\ \frac{\sqrt{-1}}{2} \left[ g^{p \bq} \left( g_{p \bk, \bq j} - g_{p \bq, \bk
j} \right) - g^{p \bm}
g_{\bm n, j} g^{n \bq} \left( g_{p \bk, \bq} - g_{p \bq, \bk} \right) \right].
\end{align*}
The result follows.
\end{proof}
\end{lemma}

\begin{lemma} \label{operatorlemma5} Given $g$ a Hermitian metric we have in
complex coordinates
\begin{align*}
\left(\del^*_{g} \del \gw \right)_{j \bk} =&\ \frac{\sqrt{-1}}{2} \left[g^{p
\bq} \left( g_{p \bk, j \bq} - g_{j \bk, p \bq} \right) + g^{p \bq}
g^{r \bs} \left( g_{p \bs, \bq} - g_{p \bq, \bs} \right) \left( g_{j \bk, r} -
g_{r \bk, j} \right) \right.\\
&\ \left. + g^{p \bq} g^{r \bs} g_{j \bq, \bs} \left( g_{p \bk, r} - g_{r \bk,
p}
\right) + g^{p \bq} g^{r \bs} g_{p \bk, \bs} \left( g_{ j \bq, r} - g_{r \bq, j}
\right) \right].
\end{align*}
\begin{proof} First of all we know that
\begin{align*}
\left(\del \gw \right)_{i j \bk} =&\ \frac{\sqrt{-1}}{2} \left(g_{j \bk, i} -
g_{i \bk, j} \right)
\end{align*}
Now, we use the general formula for $\del_{\gw}^*$ and compute
\begin{align*}
\left( \del_{g}^* \del \gw \right)_{j \bk} =&\ - g_{j \bp} g_{\bk q} \left(
\frac{\del}{\del \bar{z}^m} + \frac{1}{\bar{g}} \del_{\bm} \bar{g} \right)
\left( \del {\gw} \right)^{\bm \bp q}\\
=&\ - \frac{\sqrt{-1}}{2} \left[ g_{j \bp} g_{\bk q} \frac{\del}{\del \bar{z}^m}
\left[ g^{\bm i} g^{\bp r}
g^{\bs q} \left(g_{r \bs, i} - g_{i \bs, r} \right) \right] \right.\\
&\ \left. \qquad + g^{\bm n} g^{p \bq} g_{p \bq, \bm} \left( g_{j \bk, n} - g_{n
\bk, j} \right) \right]\\
=&\ \frac{\sqrt{-1}}{2} \left[g^{p \bq} \left( g_{p \bk, j \bq} - g_{j \bk, p
\bq} \right) \right. \\
&\ + g_{j \bp} g_{\bk q} \left( g_{r \bs, i} - g_{i \bs, r} \right) \left[
g^{\bm u} g_{u \bv, \bm} g^{i \bv} g^{\bp r} g^{\bs q} \right]\\
&\ + g_{j \bp} g_{\bk q} \left( g_{r \bs, i} - g_{i \bs, r} \right) \left[
g^{\bm i} g^{\bp u} g_{u \bv, \bm} g^{\bv r} g^{\bs q} \right]\\
&\ + g_{j \bp} g_{\bk q} \left( g_{r \bs, i} - g_{i \bs, r} \right) \left[
g^{\bm i} g^{\bp r} g^{\bs u} g_{u \bv, \bm}  g^{\bv q} \right]\\
&\ \left. - g^{\bm n} g^{p \bq} g_{p \bq, \bm} \left( g_{j \bk, n} - g_{n \bk,
j} \right) \right]\\
=&\ \frac{\sqrt{-1}}{2} \left[g^{p \bq} \left( g_{p \bk, j \bq} - g_{j \bk, p
\bq} \right) + g^{p \bq}
g^{r \bs} \left( g_{p \bs, \bq} - g_{p \bq, \bs} \right) \left( g_{j \bk, r} -
g_{r \bk, j} \right) \right.\\
&\ \left. + g^{p \bq} g^{r \bs} g_{j \bq, \bs} \left( g_{p \bk, r} - g_{r \bk,
p}
\right) + g^{p \bq} g^{r \bs} g_{p \bk, \bs} \left( g_{ j \bq, r} - g_{r \bq, j}
\right) \right].
\end{align*}
\end{proof}
\end{lemma}

\begin{lemma} \label{operatorlemma10} Given $g$ a Hermitian metric we have in
complex coordinates
\begin{gather*}
\left(\frac{\sqrt{-1}}{2} \del \delb \log \det g \right)_{j \bk} =
\frac{\sqrt{-1}}{2}  \left(g^{p
\bq} \del_j \del_{\bk} g_{p \bq} - g^{p \br}
\del_j g_{\br s} g^{s \bq}
\del_{\bk} g_{p \bq} \right)
\end{gather*}
\begin{proof}
We compute directly in coordinates
\begin{align*}
\left(\frac{\sqrt{-1}}{2} \del \delb \log \det g \right)_{j \bk} =&\
\frac{\sqrt{-1}}{2} \del_{j}
\left( g^{p \bq}
\del_{\bk} g_{p \bq} \right)\\
=&\ \frac{\sqrt{-1}}{2}  \left(g^{p \bq} \del_j \del_{\bk} g_{p \bq} - g^{p \br}
\del_j g_{\br s} g^{s \bq}
\del_{\bk} g_{p \bq} \right).
\end{align*}
\end{proof}
\end{lemma}

Also in this section we introduce canonical coordinates for $g$.  We know that
if $g$ is not K\"ahler then we cannot choose complex coordinates so that all the
first derivatives of $g$ vanish.  However, we can always ensure that a certain
symmetric part of the first derivatives vanishes.  This is made clear in the
lemma below.

\begin{lemma} \label{coords} Given a point $p \in M$, there exist coordinates
around $p$ so
that
\begin{gather*}
g_{i \bj} = \gd_{ij}
\end{gather*}
and
\begin{gather*}
\del_i g_{j \bk} + \del_j g_{i \bk} = 0.
\end{gather*}
\begin{proof} Let $\{z_i \}$ be arbitrary complex coordinate functions around
$p$ so
that $z^i(p) = 0$ for all $i$.  We briefly change our point of view and consider
the Hermitian metric $h$ associated to $g$.  The coordinate expression for $h$
takes the form
\begin{gather*}
h = h_{i j} dz^i d\bar{z}^j
\end{gather*}
where $h_{ij} = \bar{h}_{ji}$.  Without loss of generality by a rotation and
rescaling we can assume
\begin{gather} \label{coordlemma1}
h_{i j}(p) = \gd_{i j}
\end{gather}
Define new coordinates $ \{w^i \}$ by the equation
\begin{gather*}
w^i = z^i + \frac{1}{4} \sum_{j, k} \left(\frac{\del}{\del z^k} h_{i j}(p) +
\frac{\del}{\del z^j} h_{i k}(p) \right) z^j z^k
\end{gather*}
so that
\begin{gather*}
d w^i = dz^i + \frac{1}{2} \sum_{j, k} \left(\frac{\del}{\del z^k} h_{i j}(p) +
\frac{\del}{\del z^j} h_{i k }(p) \right) z^j dz^k
\end{gather*}
Note also that (\ref{coordlemma1}) still holds in these coordinates.
In these new coordinates write
\begin{gather*}
h = \til{h}_{ij} dw^i d\bar{w}^j
\end{gather*}
It is clear that
\begin{gather*}
\til{h}_{i j} = h_{ij} - \frac{1}{2} \sum_{j, k} \left(\frac{\del}{\del z^k}
h_{i j}(p) +
\frac{\del}{\del z^j}h_{i k }(p) \right) z^k + \mathcal O(z^2).
\end{gather*}
The claim follows directly by differentiating.
\end{proof}
\end{lemma}

\section{The Hermitian Hilbert Functional}
Let $(M^{2n}, g, J)$ be a complex manifold.  See section 2 for the definition of
various quantities related to the torsion.  Consider the functional
\begin{align*}
\mathbb F(g) =&\ \frac{\int_M \left[ s - \frac{1}{4} \brs{T}^2 - \frac{1}{2}
\brs{w}^2 \right] dV}{\left(\int_M dV \right)^{\frac{n-1}{n}}}.
\end{align*}

\begin{lemma} \label{functionalvar} Let $g(a)$ be a one-parameter family of
Hermitian metrics with
variation $h$.  Then
\begin{align*}
\frac{\del}{\del a} \mathbb F(g) =&\ \left( \int_M dV \right)^{\frac{1 - n}{n}}
\int_M
\left< h, - S + \frac{1}{2} Q^1 - \frac{1}{4} Q^2 - \frac{1}{2} Q^3 + Q^4
\right.\\
&\ \qquad \qquad \left. + \left( s - \frac{1}{4} \brs{T}^2 - \frac{1}{2}
\brs{w}^2 - \frac{n-1}{n} \frac{\left(\int_M s - \frac{1}{4} \brs{T}^2 -
\frac{1}{2} \brs{w}^2 \right) dV}{\int_M dV} \right) g \right> dV.
\end{align*}
Moreover, $\mathbb F$ is the \emph{unique} second-order functional which yields
$\ohat{S}$ as the leading order term in the traceless part of the variational
equation through Hermitian metrics.
\begin{proof} Combining Lemmas \ref{sfev}, \ref{q1fev} and \ref{q3fev} we see
\begin{align*}
\frac{\del}{\del a} \int_M \left[ s - \frac{1}{4} \brs{T}^2 - \frac{1}{2}
\brs{w}^2 \right] dV
=&\ \int_M \left[ \left< h, - S + \frac{1}{2} Q^1 - \frac{1}{4} Q^2 -
\frac{1}{2} Q^3 + Q^4 \right> \right.\\
&\ \qquad \left. + \tr h \left( s - \frac{1}{4} \brs{T}^2 - \frac{1}{2}
\brs{w}^2 \right) \right] dV.
\end{align*}
Likewise we compute
\begin{align*}
\frac{\del}{\del a} \left(\int_M dV \right)^{\frac{n - 1}{n}} = \frac{n-1}{n}
\int_M \tr h dV
\left( \int_M dV \right)^{-\frac{1}{n}}
\end{align*}
Combining these two calculations gives the result.  The claim of uniqueness is
also clear by inspection of the variational formulas in Lemmas \ref{sfev},
\ref{q1fev} and \ref{q3fev}.
\end{proof}
\end{lemma}

Let
\begin{gather} \label{Qdef}
Q = \frac{1}{2} Q^1 - \frac{1}{4} Q^2 - \frac{1}{2} Q^3 + Q^4
\end{gather}
 and let $K := S - Q$.  Note that
\begin{gather*}
k := \tr_g K = s - \frac{1}{4} \brs{T}^2 - \frac{1}{2} \brs{w}^2
\end{gather*}
We can rephrase the above situation in a very simple manner.  In particular
\begin{align*}
\mathbb F(g) = \int_M k dV
\end{align*}
and
\begin{align} \label{Fvar}
\frac{\del}{\del a} \mathbb F(g(a)) =&\ \int_M \left< h, - K + k g -
\frac{n-1}{n} \frac{\int_M k dV}{\int_M dV} g \right> dV
\end{align}
which is exactly analogous to the form of the gradient of the normalized Hilbert
functional.
\begin{defn} Given $(M^{2n}, g, J)$ a complex manifold we say that $g$ is
\emph{static} if $g$ is critical for $\mathbb F$.
\end{defn}

\begin{prop} Let $(M^{2n}, g, J)$ be a complex manifold with $g$ static.  Then
\begin{align*}
K - \frac{1}{n} k g = 0.
\end{align*}
Also $k$ is a constant function.  Finally, if $\mathbb F(g) \geq 0$ and $\int_M
s dV_g \leq
0$ then $g$ is K\"ahler-Einstein.
\begin{proof} The first property follows immediately by letting $h = K -
\frac{1}{n} k g$ in (\ref{Fvar}).  Next let $h = - (\gD_D k) g$ where $\gD_D$
means the Laplacian with respect to the Levi-Civita connection.  Plugging this
into (\ref{Fvar}) yields
\begin{align*}
0 =&\ - \int_M (\gD_D k) k dV = \int_M \brs{d k}^2.
\end{align*}
To see the last claim we simply note that together the hypotheses imply
\begin{align*}
0 \leq \mathbb F(g) = V^{- \frac{1}{n}} k = V^{- \frac{1}{n}} \left(s -
\frac{1}{4} \brs{T}^2 - \frac{1}{2}\brs{w}^2 \right) \leq -
\frac{1}{V^{\frac{1}{n}}} \brs{T}^2
\end{align*}
which implies $T \equiv 0$.  If the torsion of $g$ vanishes, $g$ is K\"ahler and
moreover $K$ is given by the Ricci tensor of $g$, so $g$ is K\"ahler-Einstein.
\end{proof}
\end{prop}

To emphasize, $S$ is in a sense the only natural curvature tensor associated to
a Hermitian metric which is a symmetric $(1,1)$ tensor and which is a second
order elliptic operator.  In seeking a functional which yields $S$ as the
leading term in the Euler-Lagrange equation, $\mathbb F$ above is the
\emph{only} choice.  We note that the functional $\int_M s dV_g$ was considered
in \cite{Gauduchon}. Our Lemma \ref{sfev} contains the main calculation of that
paper.  Indeed for this functional one gets automatically that critical points
are K\"ahler-Einstein if the value of the functional is nonzero.  However, there
the leading term in the Euler-Lagrange equation is $P$, which is not an elliptic
operator on Hermitian metrics.  Finally, we remark that $\mathbb F$ is
\emph{not} the Hilbert functional restricted to Hermitian metrics.  Indeed a
straightforward calculation (see also \cite{Gauduchon2} line (33)) shows that if
$r$ denotes the usual scalar curvature,
\begin{align*}
\int_M r dV_g =&\ \int_M s - \frac{1}{4} \brs{T}^2.
\end{align*}
\noindent Therefore we see that $\mathbb F = \int_M \left(r - \frac{1}{2}
\brs{w}^2 \right) dV_g$ restricted to Hermitian metrics.  This bears a certain
formal similarity to functionals related to renormalization group flows arising
in physical models \cite{OSW} \cite{Streets1}, \cite{Streets2}.

\section{Hermitian Curvature Flow}
In this section we give the definition of Hermitian curvature flow in terms of
the Chern curvature.  We then provide an equivalent definition using Hodge
operators.  In all the calculations below $Q$ is defined by (\ref{Qdef}).

\begin{prop} \label{FlowDefinition} Let $(M^{2n}, g, J)$ be a Hermitian manifold
and let
\begin{align*}
\Phi(g) := \left(S - Q \right)(g).
\end{align*}
Then $\Phi$ is a map
\begin{gather*}
\Phi(\gw) : \Re \Sym^{1,1} T^*M \to \Re \Sym^{1,1} T^* M
\end{gather*}
where $\Re \Sym^{1,1} T^* M$ are the real symmetric type $(1,1)$ tensors.
Moreover, $\Phi$ is
a nonlinear second order elliptic operator.
\begin{proof} It follows from Lemma \ref{Scalc} that
\begin{align*}
\Phi(g)_{i \bj} = - g^{k \bl} g_{i \bj, k \bl} + \mathcal O(\del g)
\end{align*}
and so $\Phi$ is a nonlinear second order elliptic operator since $g$ is
positive definite.  Also, by definition each of the tensors $Q^i$ is a real
symmetric $(1,1)$ tensor and thus $Q$ is.  It also follows from Lemma
\ref{Scalc} that $S$ is as well.  Therefore $\Phi(g)$ is a real symmetric
$(1,1)$-tensor.  The result follows.
\end{proof}
\end{prop}

\begin{defn} Given $(M^{2n}, J, g_0)$ a complex manifold with Hermitian metric
$g_0$.  We say that a one-parameter family of Hermitian metrics $g(t)$ is a
solution to \emph{Hermitian curvature flow (HCF)} with initial condition $g_0$
if
\begin{align*}
\frac{\del}{\del t} g(t) =&\ - S(g(t)) + Q(g(t))\\
g(0) =&\ g_0.
\end{align*}
\end{defn}

Next we compute a formula for HCF using Hodge operators.  Define
\begin{gather*}
\Xi = S(\cdot, J \cdot).
\end{gather*}
We will write $\Xi(g)$ using Hodge differentials.  Let
\begin{gather}
\Psi(\omega) := \del_g^* \del \omega - \del \del_g^* \omega -
\frac{\sqrt{-1}}{2} \del \delb \log \det g - 2 \sqrt{-1} \left(\delb_g^* \gw
\hook \delb \gw \right)
\end{gather}
We choose to isolate this term because as we will see in the calculations below,
this term is a real $(1,1)$ form.
\begin{prop} \label{Hodgedecomp} Given $g$ a Hermitian metric, we have
\begin{align*}
\Xi(g) =&\ \Psi(\omega) - \frac{\sqrt{-1}}{2} \left( 2 Q^4_{j \bk} + \frac{1}{2}
Q^2_{j \bk} \right).
\end{align*}
\begin{proof} Choose coordinates according to Lemma \ref{coords} so that at a
fixed
point $p \in
M$
\begin{gather*}
\del_i g_{j \bk} = - \del_j g_{i \bk}..
\end{gather*}
Using this we have that at the point $p$
\begin{gather*}
\frac{1}{2} T_{i j \bk} = \frac{1}{2} \left(\del_i g_{j \bk} - \del_j g_{i \bk}
\right) = \del_i g_{j \bk}
\end{gather*}
Next we compute a formula for $\Psi$ in coordinates.  In particular we compute a
formula for $\del \delb \log \det g$.
\begin{align*}
\frac{\sqrt{-1}}{2} \left(\del \delb \log \det g \right)_{j \bk} =&\
\frac{\sqrt{-1}}{2} \del_{j}
\left( g^{p \bq}
\del_{\bk} g_{p \bq} \right)\\
=&\ \frac{\sqrt{-1}}{2}  \left(g^{p \bq} \del_j \del_{\bk} g_{p \bq} - g^{p \br}
\del_j g_{\br s} g^{s \bq}
\del_{\bk} g_{p \bq} \right)\\
=&\ \frac{\sqrt{-1}}{2}  \left(g^{p \bq} g_{p \bq, j \bk} - \frac{1}{4} g^{p
\bq} g^{r \bs} T_{jr \bq}
T_{\bk \bs p} \right).
\end{align*}
Next we compute using Lemma \ref{operatorlemma0}
\begin{align*}
- 2 \sqrt{-1} \left(\delb_g^* \gw \hook \delb \gw \right) =&\ -2 \sqrt{-1} g^{p
\bq} \left( \delb_g^* \gw \right)_p \left(\delb \gw \right)_{\bq j \bk}\\
=&\ \frac{\sqrt{-1}}{2} g^{p \bq} g^{r \bs} T_{r p \bs} T_{\bq \bk j}.
\end{align*}
We now combine these calculations with Lemma \ref{operatorlemma1} and Lemma
\ref{operatorlemma5} to get
\begin{align*}
\Psi(\gw)_{j \bk} =&\ \frac{\sqrt{-1}}{2} \left[- g^{p \bq} g_{j \bk, p \bq}  +
\frac{1}{2} g^{p \bq} g^{r \bs} \left[ 2 T_{\bq \bs p} T_{r j \bk} + 2 T_{p r
\bq} T_{\bs \bk j} \right. \right.\\
&\ \left. \left. + T_{\bs
\bq j} T_{r p \bk} + T_{\bs \bk p} T_{r j \bq} - T_{j r \bq} T_{\bk \bs p} +
\frac{1}{2}
T_{j r
\bq} T_{\bk \bs p} \right] \right]\\
=&\ \frac{\sqrt{-1}}{2}  \left[- g^{p \bq} g_{j \bk, p \bq} \right.\\
&\ \left.+ \frac{1}{2} g^{p
\bq} g^{r \bs} \left[ 2 T_{\bq
\bs p} T_{r j \bk} + 2 T_{p r \bq} T_{\bs \bk j} + T_{\bs \bq j} T_{r p \bk} +
\frac{1}{2} T_{j r \bq} T_{\bk
\bs p} \right] \right].
\end{align*}
Likewise we have from Lemma \ref{Scalc}
\begin{align*}
\Xi(g)_{j \bk} =&\ \frac{\sqrt{-1}}{2} \left[- g^{p \bq} g_{j \bk, p \bq} +
\frac{1}{4} g^{p \bq} g^{r
\bs} T_{j p \bs} T_{\bk \bq r}\right].
\end{align*}
The result follows by combining these calculations.
\end{proof}
\end{prop}

\begin{cor} \label{HCFCC} The HCF equation is equivalent to
\begin{gather*}
\begin{split}
\frac{\del}{\del t} \omega =&\ - \Psi(\omega) + \frac{\sqrt{-1}}{2} \left(
\frac{1}{2} Q^1 + \frac{1}{4} Q^2 - \frac{1}{2} Q^3 + 3 Q^4 \right).
\end{split}
\end{gather*}
\begin{proof} This follows immediately from the definition of HCF and
Proposition \ref{Hodgedecomp}.
\end{proof}
\end{cor}

\section{Short-Time Existence}

\begin{prop} \label{ste} Given $(M^{2n}, J, g_0)$ a compact complex manifold,
there exists a
unique solution to HCF with initial condition $g_0$ on $[0, \ge)$ for some $\ge
> 0$.
\begin{proof} Since the operator $\Phi(g)$ is strictly elliptic by Proposition
\ref{FlowDefinition} the HCF equation is strictly parabolic, and thus short-time
existence and uniqueness follows from standard theory.
\end{proof}
\end{prop}

\begin{prop} Given $(M^{2n}, J, g_0)$ a compact complex manifold with K\"ahler
metric $g_0$, let $g(s)$ denote the solution to HCF with initial condition
$g_0$, which exists on $[0, T)$.  Then for all $t \in [0, T)$, $g(t)$ is
K\"ahler, and is a solution to K\"ahler Ricci flow.
\begin{proof} Let $\til{g}(t)$ be the solution to K\"ahler Ricci flow with
initial condition $g_0$.  Ricci flow preserves the K\"ahler condition, thus
$\til{g}(t)$ is K\"ahler for all time, hence $\til{T} = d \til{\omega} = 0$ and
$\til{S} = \til{\Rc}$.  Together this implies that $\til{g}(t)$ satisfies
\begin{align*}
\frac{\del}{\del t} \til{g}(t) =&\ - \til{\Rc}\\
=&\ - \til{S} + \til{Q}
\end{align*}
Thus $\til{g}(t)$ is a solution to HCF with initial condition $g_0$.  Since
solutions to HCF are unique, it follows that $\til{g}(t) = g(t)$ for all time
and hence $g(t)$ is K\"ahler for all time and solves K\"ahler-Ricci flow.
\end{proof}
\end{prop}

\section{Evolution Equations}
\begin{lemma} \label{Omegaev} For a solution to HCF we have
\begin{align*}
\dt \Omega_{i \bj k \bl} =&\ \gD \Omega_{i \bj k \bl} + g^{m \bn} \left(T_{\bn
\bj}^{\bp} \N_m \Omega_{i {\bp} k \bl}  + T_{m i}^p \N_{\bj} \Omega_{p \bn k
\bl} \right)\\
&\ + g^{m \bn} \left(\Omega_{i \bj m}^p \Omega_{p \bn k \bl} + \Omega_{m \bn
\bj}^{\bp} \Omega_{i \bp k \bl} + \Omega_{m \bj k}^p \Omega_{i \bn p \bl} +
\Omega_{m \bj \bl}^{\bp} \Omega_{i \bn k \bp} \right)\\
&\ - \Omega_{i \bj k}^m \left(S_{m \bl} - Q_{m \bl} \right) - \N_{\bj} \N_i Q_{k
\bl}
\end{align*}
\begin{proof} First consider the term $Q$ in the evolution of $g$.  Using Lemma
\ref{curvaturevariation} we see that this contributes
\begin{align*}
\Omega_{i \bj k}^m Q_{m \bl} - \N_{\bj} \N_i Q_{k \bl}.
\end{align*}
to the evolution of $\Omega_{i \bj k \bl}$.  Next we consider the contribution
of the term $-S$ in the evolution of $g$.  Using Lemma \ref{curvaturevariation}
we see that the evolution $\dt g = - S$ yields
\begin{align*}
\frac{\del}{\del t} \Omega_{i \bj k \bl} =&\ - \Omega_{i \bj k}^m S_{m \bl} +
\N_{\bj} \N_i S_{k \bl}.
\end{align*}
Now we must apply the second Bianchi identity.  We have
\begin{align*}
\N_{\bj} \left( \N_i S_{k \bl} \right) =&\ \N_{\bj} g^{m \bn} \left( \N_i
\Omega_{m \bn k \bl} \right)\\
=&\ \N_{\bj} g^{m \bn} \left( \N_m \Omega_{i \bn k \bl} + T_{m i}^p \Omega_{p
\bn k \bl} \right)\\
=&\ g^{m \bn} \left(\N_{\bj} \N_m \Omega_{i \bn k \bl} + \N_{\bj} T_{m i}^p
\Omega_{p \bn k \bl} + T_{m i}^p \N_{\bj} \Omega_{p \bn k \bl} \right)
\end{align*}
Next we commute covariant derivatives to get
\begin{align*}
g^{m \bn} \N_{\bj} \N_m \Omega_{i \bn k \bl} =&\ g^{m \bn} \left(\N_m \N_{\bj}
\Omega_{i \bn k \bl} + \Omega_{m \bj i}^p \Omega_{p \bn k \bl} + \Omega_{m \bj
\bn}^{\bp} \Omega_{i \bp k \bl} \right.\\
&\ \left. \qquad + \Omega_{m \bj k}^p \Omega_{i \bn p \bl} + \Omega_{m \bj
\bl}^{\bp} \Omega_{i \bn k \bp} \right).
\end{align*}
Finally, we apply the Bianchi identity again to get
\begin{align*}
g^{m \bn} \N_m \N_{\bj} \Omega_{i \bn k \bl} =&\ g^{m \bn} \N_m \left(\N_{\bn}
\Omega_{i \bj k \bl} + T_{\bn \bj}^{\bp} \Omega_{i \bp k \bl} \right)\\
=&\ \gD \Omega_{i \bj k \bl} + g^{m \bn} \left( \N_m T_{\bn \bj}^{\bp} \Omega_{i
\bp k \bl} + T_{\bn \bj}^{\bp} \N_m \Omega_{i p k \bl} \right)
\end{align*}
Combining these calculations yields
\begin{align*}
\dt \Omega_{i \bj k \bl} =&\ \gD \Omega_{i \bj k \bl}\\
&\ + g^{m \bn} \left(\N_m T_{\bn \bj}^{\bp} \Omega_{i \bp k \bl} + T_{\bn
\bj}^{\bp} \N_m \Omega_{i {\bp} k \bl} + \N_{\bj} T_{m i}^p \Omega_{p \bn k \bl}
+ T_{m i}^p \N_{\bj} \Omega_{p \bn k \bl} \right)\\
&\ + g^{m \bn} \left(\Omega_{m \bj i}^p \Omega_{p \bn k \bl} + \Omega_{m \bj
\bn}^{\bp} \Omega_{i \bp k \bl} + \Omega_{m \bj k}^p \Omega_{i \bn p \bl} +
\Omega_{m \bj \bl}^{\bp} \Omega_{i \bn k \bp} \right) - \Omega_{i \bj k}^m S_{m
\bl}
\end{align*}
Now we can apply the Bianchi identity to the terms
\begin{align*}
\Omega_{m \bj \bn}^{\bp} =&\ \Omega_{m \bn \bj}^{\bp}
+ \N_m T_{\bj \bn}^{\bp}
\end{align*}
and
\begin{align*}
\Omega_{m \bj i}^p =&\ \Omega_{i \bj m}^p + \N_{\bj} T_{i m}^p.
\end{align*}
Plugging these two in yields
\begin{align*}
\dt \Omega_{i \bj k \bl} =&\ \gD \Omega_{i \bj k \bl} + g^{m \bn} \left(T_{\bn
\bj}^{\bp} \N_m \Omega_{i {\bp} k \bl}  + T_{m i}^p \N_{\bj} \Omega_{p \bn k
\bl} \right)\\
&\ + g^{m \bn} \left(\Omega_{i \bj m}^p \Omega_{p \bn k \bl} + \Omega_{m \bn
\bj}^{\bp} \Omega_{i \bp k \bl} + \Omega_{m \bj k}^p \Omega_{i \bn p \bl} +
\Omega_{m \bj \bl}^{\bp} \Omega_{i \bn k \bp} \right)\\
&\ - \Omega_{i \bj k}^m S_{m \bl}.
\end{align*}
Combining this with the above terms gives the result.
\end{proof}
\end{lemma}

\begin{lemma} \label{Tev} For a solution to HCF we have
\begin{align*}
\dt T_{i j \bk} =&\ \gD T_{i j \bk} + g^{m \bn} \left[T_{j i}^p \N_{\bn} T_{m p
\bk} + \N_{\bn} T_{m j}^p T_{i p
\bk} + T_{m j}^p \N_{\bn} T_{i p \bk} \right.\\
&\ \left. \qquad \quad + \N_{\bn} T_{i m}^p T_{j p \bk} + T_{i m}^p \N_{\bn}
T_{j p \bk} \right]\\
&\ + g^{m \bn} \left[\Omega_{\bn j m}^p T_{i p \bk} + \Omega_{\bn j \bk}^{\bp}
T_{i m \bp} - \Omega_{\bn i m}^p T_{j p \bk} \right.\\
&\ \qquad \quad \left. - \Omega_{\bn i \bk}^{\bp} T_{j m \bp} - \Omega_{p \bn m
\bk} T_{j i}^p \right] - T_{i j}^p \left(S_{p \bk} - Q_{p \bk} \right)\\
&\ + \N_i Q_{j \bk} - \N_j Q_{i \bk}.
\end{align*}
\begin{proof} First we compute the contribution from the term $Q$ in the
evolution of $g$.  In particular using Lemma \ref{torsionvariation} this yields
\begin{align*}
\N_i Q_{j \bk} - \N_j Q_{i \bk} + T_{i j}^p Q_{p \bk}.
\end{align*}
Next we focus on the term $-S$.  Applying Lemma \ref{torsionvariation} we have.
\begin{align*}
\N_j S_{i \bk} - \N_i S_{j \bk} - T_{i j}^p S_{p \bk}.
\end{align*}
Now we rewrite using the Bianchi identity
\begin{align*}
\N_j S_{i \bk} =&\ g^{m \bn} \N_j \Omega_{m \bn i \bk}\\
=&\ g^{m \bn} \N_j \left( \Omega_{i \bn m \bk} + \N_{\bn} T_{i m \bk} \right)
\end{align*}
and
\begin{align*}
\N_i S_{j \bk} =&\ g^{m \bn} \N_i \left(\Omega_{j \bn m \bk} + \N_{\bn} T_{j m
\bk} \right).
\end{align*}
Combining these yields
\begin{align*}
\dt T_{i j \bk} =&\ g^{m \bn} \left( \N_j \Omega_{i \bn m \bk} - \N_i \Omega_{j
\bn m \bk} + \N_j \N_{\bn} T_{i m \bk} - \N_i \N_{\bn} T_{j m \bk} \right)\\
&\ - T_{i j}^p S_{p \bk}.
\end{align*}
Applying the Bianchi identity again yields
\begin{align*}
g^{m \bn} \left( \N_j \Omega_{i \bn m \bk} - \N_i \Omega_{j \bn m \bk} \right)
=&\ - g^{m \bn} T_{j i}^p \Omega_{p \bn m \bk}.
\end{align*}
Also, we commute derivatives
\begin{align*}
\N_j \N_{\bn} T_{i m \bk} =&\ \N_{\bn} \N_j T_{i m \bk} + \Omega_{\bn j i}^p
T_{p m \bk} + \Omega_{\bn j m}^p T_{i p \bk} + \Omega_{\bn j \bk}^{\bp} T_{i m
\bp}
\end{align*}
and
\begin{align*}
\N_i \N_{\bn} T_{j m \bk} =&\ \N_{\bn} \N_i T_{j m \bk} + \Omega_{\bn i j}^p
T_{p m \bk} + \Omega_{\bn i m}^p T_{j p \bk} + \Omega_{\bn i \bk}^{\bp} T_{j m
\bp}.
\end{align*}
Finally, using Lemma \ref{Tformula} we see
\begin{align*}
g^{m \bn} \N_{\bn} \left( \N_j T_{i m \bk} - \N_i T_{j m \bk} \right) =&\ g^{m
\bn} \N_{\bn} \left( \N_m T_{i j \bk} + T_{j i}^p T_{m p \bk} + T_{m j}^p T_{i p
\bk} + T_{i m}^p T_{j p \bk} \right)\\
=&\ \gD T_{i j \bk} + g^{m \bn} \N_{\bn} \left(T_{j i}^p T_{m p \bk} + T_{m j}^p
T_{i p \bk} + T_{i m}^p T_{j p \bk} \right).
\end{align*}
Combining these calculations yields
\begin{align*}
\dt T_{i j \bk} =&\ \gD T_{i j \bk} + g^{m \bn} \N_{\bn} \left(T_{j i}^p T_{m p
\bk} + T_{m j}^p T_{i p \bk} + T_{i m}^p T_{j p \bk} \right)\\
&\ + g^{m \bn} \left[\Omega_{\bn j i}^p T_{p m \bk} + \Omega_{\bn j m}^p T_{i p
\bk} + \Omega_{\bn j \bk}^{\bp} T_{i m \bp} \right.\\
&\ \left. - \Omega_{\bn i j}^p T_{p m \bk} - \Omega_{\bn i m}^p T_{j p \bk} -
\Omega_{\bn i \bk}^{\bp} T_{j m \bp} - \Omega_{p \bn m \bk} T_{j i}^p \right] -
T_{i j}^p S_{p \bk}.
\end{align*}
Using the Bianchi identity we can simplify
\begin{align*}
g^{m \bn} \left[ \Omega_{\bn j i}^p T_{p m \bk} - \Omega_{\bn i j}^p T_{p m \bk}
\right] =&\ g^{m \bn} \left[ \Omega_{\bn j i}^p + \Omega_{i \bn j}^p \right]
T_{p m \bk}\\
=&\ g^{m \bn} \N_{\bn} T_{j i}^p T_{p m \bk}.
\end{align*}
Plugging this in yields
\begin{align*}
\dt T_{i j \bk} =&\ \gD T_{i j \bk} + g^{m \bn} \left[T_{j i}^p \N_{\bn} T_{m p
\bk} + \N_{\bn} T_{m j}^p T_{i p
\bk} + T_{m j}^p \N_{\bn} T_{i p \bk} \right.\\
&\ \left. \qquad \quad + \N_{\bn} T_{i m}^p T_{j p \bk} + T_{i m}^p \N_{\bn}
T_{j p \bk} \right]\\
&\ + g^{m \bn} \left[\Omega_{\bn j m}^p T_{i p \bk} + \Omega_{\bn j \bk}^{\bp}
T_{i m \bp} - \Omega_{\bn i m}^p T_{j p \bk} \right.\\
&\ \qquad \quad \left. - \Omega_{\bn i \bk}^{\bp} T_{j m \bp} - \Omega_{p \bn m
\bk} T_{j i}^p \right] - T_{i j}^p S_{p \bk}.
\end{align*}
Combining this with the terms from $Q$ gives the result.
\end{proof}
\end{lemma}

\section{Higher Derivative Estimates}
In this section we will prove derivative estimates for HCF.  It will be most
convenient to phrase these results in terms of the curvature of the Chern
connection.  All of the calculations below will be done in canonical coordinates
at a fixed point.  In particular, in these coordinates any first derivative of
$g$ can be expressed in terms of the torsion $T$, and any second derivative can
be expressed in terms of a sum of curvature and torsion.

\begin{lemma} \label{CurDerEv} Given $(M^{2n}, g(t), J)$ a solution to HCF we
have
\begin{align*}
\dt \N^k \Omega =&\ \gD \N^k \Omega + \sum_{j = 0}^k \N^j T * \N^{k + 1 - j}
\Omega + \sum_{j=0}^k \N^j
\Omega *
\N^{k-j} \Omega\\
&\ + \sum_{j=0}^{k} \sum_{l=0}^j \N^l T * \N^{j-l} T * \N^{k-j} \Omega.
\end{align*}
\begin{proof} The case $k = 0$ is covered by Lemma \ref{Omegaev}.  We directly
compute
\begin{align*}
\dt \N^k \Omega =&\ \dt \left( \del + \gG \right) * \left( \del + \gG \right) *
\dots * \left(\del + \gG \right) \Omega\\
=&\ \N^k \left( \dt \Omega \right) + \left( \dt \gG \right) * \left( \del + \gG
\right) \dots * \left(\del + \gG \right) \Omega\\
&\ + \left( \del + \gG \right) * \left( \dt \gG \right) * \dots * \left( \del +
\gG \right) \Omega + \dots\\
&\ + \left( \del + \gG \right) * \dots * \left(\del + \gG \right) * \left(\dt
\gG \right) \Omega
\end{align*}
We apply Lemma \ref{connectionvariation} to see that
\begin{align*}
\dt \gG =&\ \N \left(\Omega + T^{*2} \right) = \N \Omega + T * \N T.
\end{align*}
Plugging this in yields
\begin{align*}
\dt \N^k \Omega =&\ \N^k \left( \gD \Omega + T * \N \Omega + \Omega^{*2} +
\Omega * T^{*2}
\right)\\
&\ + \sum_{j=0}^{k-1} \N^j \left(\N \Omega  + T * \N T \right) * \N^{k-1-j}
\Omega\\
=&\ \gD \N^k \Omega + \sum_{j = 0}^k \N^j T * \N^{k + 1 - j} \Omega +
\sum_{j=0}^k \N^j \Omega * \N^{k-j}
\Omega\\
&\ + \sum_{j=0}^{k} \sum_{l=0}^j \N^l T * \N^{j-l} T * \N^{k-j} \Omega.
\end{align*}
\end{proof}
\end{lemma}

\begin{lemma} \label{torsderev} Given $(M^{2n}, g(t), J)$ a solution to HCF we
have
\begin{align*}
\dt \N^k T =&\ \gD \N^k T + \sum_{j = 0}^{k+1} \N^j T * \N^{k+1 - j} T + \sum_{j
= 0}^{k} \N^j T * \N^{k - j} \Omega\\
&\ + \sum_{j = 0}^{k-1} \sum_{l = 0}^j \N^l T * \N^{j - l + 1} T * \N^{k - 1 -
j} T.
\end{align*}
\begin{proof} The case $k = 0$ is covered by Lemma \ref{Tev}.  We directly
compute
\begin{align*}
\dt \N^k T =&\ \dt \left( \del + \gG \right) * \dots * \left(\del + \gG \right)
T\\
=&\ \N^k \left(\dt T \right) + \left( \dt \gG \right) * \left( \del + \gG
\right) \dots * \left(\del + \gG \right) T\\
&\ + \left( \del + \gG \right) * \left( \dt \gG \right) * \dots * \left( \del +
\gG \right) T + \dots\\
&\ + \left( \del + \gG \right) * \dots * \left(\del + \gG \right) * \left(\dt
\gG \right) T.
\end{align*}
Again we apply Lemma \ref{connectionvariation} to yield
\begin{align*}
=&\ \N^k \left( \gD T + \N T * T + \Omega * T \right)\\
&\ + \sum_{j = 0}^{k-1} \N^j \left(\N \Omega + T * \N T \right) * \N^{k - 1 - j}
T\\
=&\ \gD \N^k T + \sum_{j = 0}^{k+1} \N^j T * \N^{k+1 - j} T + \sum_{j = 0}^{k}
\N^j T * \N^{k - j} \Omega\\
&\ + \sum_{j = 0}^{k-1} \sum_{l = 0}^j \N^l T * \N^{j - l + 1} T * \N^{k - 1 -
j} T.
\end{align*}
\end{proof}
\end{lemma}

\begin{thm} \label{derivthm1} Let $(M^{2n}, g(t), J)$ be a solution to
HCF for which the maximum principle holds.  Then for each $\ga > 0$ and
every $m \in \mathbb N$ there exists a constant $C_m$ depending only on $m, n$
and $\max \{\ga, 1 \}$ such that if
\begin{gather} \label{derivloc10}
\begin{split}
\brs{\Omega}_{C^0(g_t)} \leq&\ K,\\
\brs{\N T}_{C^0(g_t)} \leq&\ K,\\
\brs{T}^2_{C^0(g_t)} \leq&\ K
\end{split}
\end{gather}
for all $x \in M$ and $t \in \left[0, \frac{\ga}{K} \right]$, then
\begin{gather} \label{derivloc20}
\begin{split}
\brs{\N^m \Omega}_{C^0(g_t)} \leq&\ \frac{C_m K}{t^{m/2}},\\
\brs{\N^{m+1} T}_{C^0(g_t)} \leq&\ \frac{C_m K}{t^{m/2}}
\end{split}
\end{gather}
for all $x \in M$ and $t \in \left(0, \frac{\ga}{K} \right]$.
\begin{proof} Our proof is by induction on $m$.  First consider $m = 1$.  The
following evolution equation for $\brs{\Omega}^2$ follows from Lemma
\ref{CurDerEv}:
\begin{gather} \label{derivloc22}
\dt \brs{\Omega}^2 = \gD \brs{\Omega}^2 - 2 \brs{\N \Omega}^2 + T * \N \Omega *
\Omega +
\Omega^{*3} + T^{*2} * \Omega^{*2}.
\end{gather}
Also from Lemma \ref{CurDerEv} we conclude
\begin{gather} \label{derivloc21}
\begin{split}
\dt \brs{\N \Omega}^2 =&\ \gD \brs{\N \Omega}^2 - 2 \brs{\N^2 \Omega}^2  + T *
\N^2 \Omega * \N \Omega\\
&\ + \Omega * \N \Omega^{*2} + \N T * \N \Omega^{*2} + T^{*2} * \N \Omega^{*2}\\
&\ + T * \N T * \Omega * \N \Omega.
\end{split}
\end{gather}
Now, we aim to use the term $ -2 \brs{\N \Omega}^2$ in the evolution of
$\brs{\Omega}^2$ to control the evolution of $\brs{\N \Omega}^2$.  Consider the
function
\begin{gather*}
F(x,t) := t \brs{\N \Omega}^2 + \gb \brs{\Omega}^2
\end{gather*}
where $\gb$ is a constant to be chosen below.  Putting together
(\ref{derivloc22})
and (\ref{derivloc21}) gives
\begin{gather} \label{derivloc30}
\begin{split}
\dt F \leq&\ \gD F - 2 t \brs{\N^2 \Omega}^2 + \left(1 + c_1 t \brs{\Omega} - 2
\gb
\right) \brs{\N \Omega}^2\\
&\ + t \left( T * \N^2 \Omega * \N \Omega + \Omega * \N \Omega^{*2} + \N T * \N
\Omega^{*2} \right.\\
&\ \qquad \left. + T^{*2} * \N \Omega^{*2} + T * \N T * \Omega * \N \Omega  + T
* \Omega^{*2} * \N \Omega \right)\\
&\ + c_2 \gb \left( \brs{T}^2 \brs{\Omega}^2 + \brs{\Omega}^3 \right) + T * \N
\Omega * \Omega
\end{split}
\end{gather}
where all the $c_i$ are universal constants depending only on dimension.  We
must estimate the different terms in (\ref{derivloc30}).  First of all we use
the
Cauchy-Schwarz inequality and the assumption (\ref{derivloc10}) to conclude
\begin{gather} \label{derivloc31}
\begin{split}
t T * \N^2 \Omega * \N \Omega \leq&\ t c_3 \left( \brs{\N^2 \Omega} \right)
\left(
\brs{T} \brs{\N \Omega} \right)\\
\leq&\ t c_3 \left( \frac{\brs{\N^2 \Omega}^2}{2 c_3} + \frac{c_3 \brs{T}^2
\brs{\N \Omega}^2}{2} \right)\\
\leq&\ \frac{t}{2} \brs{\N^2 \Omega}^2 + c_4 t K \brs{\N \Omega}^2.
\end{split}
\end{gather}
Similarly we simplify
\begin{gather} \label{derivloc32}
\begin{split}
t T * \Omega^{*2} * \N \Omega \leq&\ t c_5 \left(\brs{\Omega}^2 \right)
\left(\brs{T} \brs{\N \Omega} \right)\\
\leq&\ c_6 t K^4 + c_6 t K \brs{\N \Omega}^2
\end{split}
\end{gather}
and also
\begin{gather} \label{derivloc33}
T * \Omega * \N \Omega \leq \frac{c_7}{\gb} K^3 + \gb \brs{\N \Omega}^2.
\end{gather}
The remaining terms are estimated in an obviously analogous fashion.  Plugging
these estimates into (\ref{derivloc30}) gives
that
for $t \in \left[0, \frac{\ga}{K} \right]$,
\begin{gather*}
\dt F \leq \gD F + \left(1 + c_8 \ga - \gb \right) \brs{\N \Omega}^2 + c_9
\left(\ga + \gb \right) K^3.
\end{gather*}
Choose $\gb \geq \frac{1 + c_8 \ga}{2}$ and note that $\gb$ depends only on the
dimension and $\max \{\ga, 1 \}$.  Then we have that for $t \in \left[0,
\frac{\ga}{K} \right]$,
\begin{gather*}
\dt F \leq c_{10} \gb K^3
\end{gather*}
Using that $F(0) \leq \gb K^2$ and applying the maximum principle gives
\begin{gather*}
\sup_{x \in M} F(x,t) \leq \gb K^2 + c_{10} \gb K^3 t \leq \left(1 + c_{10} \ga
\right) \gb K^2 \leq C_1^2 K^2
\end{gather*}
where again $C_1$ depends only on $n$ and $\max \{\ga, 1 \}$.  Thus
\begin{gather*}
\brs{\N \Omega} \leq \sqrt{ \frac{F}{t}} \leq \frac{C_1 K}{t^{1/2}}
\end{gather*}
for all $x \in M$ and $t \in \left(0, \frac{\ga}{K} \right]$.  To get the
estimate for $\N^2 T$ one computes the evolution of a function
\begin{gather*}
F(x,t) := t \left( \brs{\N \Omega}^2 + \brs{\N^2 T}^2 \right) + \gb
\brs{\Omega}^2
\end{gather*}
and argues as above, the bounds being entirely analogous.  This completes the
case $m = 1$.  

For the induction step we first conclude from Lemma
\ref{CurDerEv} the evolution equation
\begin{align*}
\dt \brs{\N^k \Omega}^2 =&\ \gD \brs{\N^k \Omega}^2 - 2 \brs{\N^{k+1} \Omega}^2
+ \sum_{j = 0}^{k} \N^j T * \N^{k+1 - j} \Omega * \N^{k} \Omega\\
&\ + \sum_{j=0}^k \N^j \Omega * \N^{k-j} \Omega * \N^k \Omega\\
&\ + \sum_{j = 0}^k \sum_{l = 0}^j \N^l T * \N^{j - l} T * \N^{k - j} \Omega *
\N^k \Omega.
\end{align*}
We address the first sum in the above equation.  We first make the bound
\begin{gather} \label{derivloc40}
\begin{split}
T * \N^{k+1} \Omega * \N^k \Omega \leq&\ c \brs{T} \brs{\N^{k+1} \Omega}
\brs{\N^k
\Omega}\\
\leq&\ \ge \brs{\N^{k+1} \Omega}^2 + C(\ge) \brs{T}^2 \brs{\N^k \Omega}^2\\
\leq&\ \ge \brs{\N^{k+1} \Omega}^2 + C(\ge) K \brs{\N^k \Omega}^2.
\end{split}
\end{gather}
Also we have
\begin{gather}
\begin{split}
\N T * \N^k \Omega * \N^k \Omega \leq&\ c \brs{\N T} \brs{\N^k \Omega}^2\\
\leq&\ K \brs{\N^k \Omega}^2.
\end{split}
\end{gather}
For the rest of the summand we bound for $j > 0$
\begin{gather}
\begin{split}
\N^j T * \N^{k + 1 - j} \Omega * \N^k \Omega \leq&\ c \brs{\N^j T} \brs{\N^{k +
1 - j} \Omega} \brs{\N^k \Omega}\\
\leq&\ c \frac{K}{t^{(j - 1)/2}} \frac{K}{t^{(k + 1 - j)/2}} \brs{\N^k \Omega}\\
\leq&\ c K \brs{\N^k \Omega}^2 + c \frac{K^3}{t^k}.
\end{split}
\end{gather}
A similar calculation yields a bound
\begin{gather}
\begin{split}
\N^j \Omega * \N^{k-j} \Omega * \N^k \Omega \leq&\ c \brs{\N^j \Omega}
\brs{\N^{k-j} \Omega}
\brs{\N^k \Omega}\\
\leq&\ c \frac{K}{t^{j/2}} \frac{K}{t^{(k-j)/2}} \brs{\N^k \Omega}\\
\leq&\ c K \brs{\N^k \Omega}^2 + c \frac{K^3}{t^k}.
\end{split}
\end{gather}
Next we bound using the inequality $K \leq \frac{C}{t}$
\begin{gather}\label{derivloc45}
\begin{split}
\N^l T * \N^{j - l} T * \N^{k - j} \Omega * \N^k \Omega \leq&\ c \brs{\N^l T}
\brs{\N^{j-l} T} \brs{\N^{k - j} \Omega} \brs{\N^k \Omega}\\
\leq&\ \frac{K}{t^{(l-1)/2}} \frac{K}{t^{(j - l - 1)/2}} \frac{K}{t^{(k-j)/2}}
\brs{\N^k \Omega}\\
\leq&\ \frac{K^3}{t^{(k-2)/2}} \brs{\N^k \Omega}\\
\leq&\ \frac{K^2}{t^{(k-1)/2}} \brs{\N^k \Omega}\\
\leq&\ c K \brs{\N^k \Omega}^2 + c \frac{K^3}{t^k}.
\end{split}
\end{gather}
Using (\ref{derivloc40}) - (\ref{derivloc45}) we conclude
\begin{gather*}
\dt \brs{\N^k \Omega}^2 \leq \gD \brs{\N^k \Omega}^2 - \brs{\N^{k+1} \Omega}^2 +
C K
\left( \brs{\N^k \Omega}^2 + \frac{K^2}{t^k} \right).
\end{gather*}
Furthermore, using completely analogous bounds one can conclude
\begin{gather*}
\dt \brs{\N^{k+1} T}^2 \leq \gD \brs{\N^{k+1} T}^2 - \brs{\N^{k+2} T}^2 +
\frac{1}{2} \brs{\N^{k+1} \Omega}^2 + C K \left( \brs{\N^{k+1} T}^2 +
\frac{K^2}{t^k} \right).
\end{gather*}
The extra term $\frac{1}{2} \brs{\N^{k + 1} \Omega}^2$ arises from the term $T *
\N^{k+1} \Omega * \N^{k+1} T$.  Together these yield, if we set $H_k = \brs{\N^k
\Omega}^2 + \brs{\N^{k+1} \Omega}^2$,
\begin{align*}
\dt H_k \leq&\ \gD H_k - \frac{1}{2} H_{k+1} + C K H_k + \frac{K^3}{t^k}.
\end{align*}
This bound is sufficient to carry out the inductive step analogously to the step
$k = 1$.  The details of this construction are found in \cite{Chow} page
229-230.
\end{proof}
\end{thm}

\begin{cor} \label{STECor} There exists a constant $c = c(n)$ such that given
$(M^{2n}, g, J)$ a complex manifold with Hermitian metric $g$, the solution to
HCF with initial condition $g$ exists for $t \in \left[0, \frac{c(n)}{\max \{
\brs{\Omega}_{C^0}, \brs{\N T}_{C^0}, \brs{T}_{C_0}^2 \}} \right]$.  Moreover
the solution exists on a maximal time interval $[0, \tau)$, and if $\tau <
\infty$ then
\begin{align*}
\limsup_{t \to \tau} \max \{ \brs{\Omega}_{C^0(g_t)}, \brs{\N T}_{C^0(g_t)},
\brs{T}^2_{C^0(g_t)} \} = \infty.
\end{align*}
\begin{proof} This argument is standard.  Using the evolution equations for $T$,
$\N T$
and $\Omega$ it is easy to prove a ``doubling-time'' estimate for these
quantities on the interval stated using the maximum principle.  Once this is in
place, the derivative
estimates follow from Theorem \ref{derivthm1}.  These yield bounds on the
curvature and torsion and all covariant derivatives on the stated interval,
which can be integrated in time to show smooth existence of the flow on that
interval.

Finally, if one has that the curvature, torsion and first covariant derivative
of torsion are bounded up to a time $\tau < \infty$, one concludes from Theorem
\ref{derivthm1} uniform bounds on the derivatives of curvature and torsion on
$[0, \tau]$.  These bounds can be integrated in time to get $C^k$ bounds on the
metric on this whole time interval, yielding smooth existence up to this time.
\end{proof}
\end{cor}

\noindent Note now that Theorem \ref{STEThm} is a consequence of Proposition
\ref{ste}, Theorem \ref{derivthm1} and Corollary \ref{STECor}.

\section{Stability}
In this section we prove dynamic stability of HCF near a K\"ahler-Einstein
metric with negative or zero first Chern class\footnote{In an earlier version of
this paper we claimed stability for positive first Chern class as well.
However, such a stability result needs a further condition on initial metric,
since for instance
the K\"ahler class under volume-normalized Ricci flow satisfies $[\omega_t] =
c_1(M) + e^{t} \left( [\omega_0] - c_1(M) \right)$. Thus even nearby
K\"ahler metrics will not converge under the resulting K\"ahler Ricci flow.
This is reflected in the fact that K\"ahler-Einstein metrics with positive first
Chern class are not linearly stable, which is what causes our proof to fail.  We
expect that for a given K\"ahler-Einstein metric $g$ with $c_1(M) >
0$ there exist certain harmonic (1,1)-forms $h$ such that the solution to HCF
with
initial condition $g + h$ converges modulo automorphisms to a K\"ahler-Einstein
metric.  This will be the subject of future work.}.  By examining the linearized
deformation equation we know that
K\"ahler-Einstein metrics are rigid in case where $c_1(M) < 0$.  In the case
$c_1(M) = 0$, there can be nontrivial deformation of
K\"ahler-Einstein metrics due to variation of K\'ahler class.  There is a
general technique for
dealing with stability of evolution equations around integrable stationary
points \cite{Tian}, \cite{Sesum}, \cite{Simon}.  Given the discussion above, our
problem falls squarely into the realm of these techniques, and so we adopt them.
 We note that since the $c_1(M) < 0$ case is rigid, there may be an easier proof
for this case, but in the interest of covering the most cases possible with a
single proof we choose the more general technique.

Consider the volume-normalized HCF equation
\begin{align*}
\dt g =&\ - S + Q + \frac{1}{n} \left( \int_M \tr_g \left(S - Q \right) dV
\right) g\\
=&:\ - \mathcal F(g).
\end{align*}
We compute the linearization of $\mathcal F$ around a K\"ahler-Einstein
metric.
Since the tensor $\mathcal F(g)$ is only defined for Hermitian metrics we
obviously compute the
variation of $\mathcal F(g)$ through a family of Hermitian metrics.

\begin{prop} \label{LDE} Let $(M^{2n}, J)$ be a complex manifold and suppose
$g(a)$ is a
one-parameter family of unit volume Hermitian metrics compatible with $J$ with
\begin{gather*}
\frac{\del}{\del a} g(a) |_{a = 0} = h.
\end{gather*}
Moreover suppose $g(0)$ is K\"ahler-Einstein.  Then
\begin{align*}
\frac{\del}{\del a} \mathcal F(g) = \N^* \N h - \overset{\circ}{R}(h)
\end{align*}
where $\overset{\circ}{R}(h)_{k \bl} = h^{i \bj} R_{k \bj i \bl}$.
\begin{proof}  Choose complex coordinates which are normal for $g(0)$ at a point
$p \in
M$.  First we note that
\begin{gather*}
\begin{split}
\left. \frac{\del}{\del a} T(a) * T(a) \right|_{a = 0} =&\ h * T(0) * T(0) +
\left( \frac{\del}{\del a} T(a) \right) * T(0)\\
=&\ 0
\end{split}
\end{gather*}
since the metric $g(0)$ is K\"ahler and hence torsion-free.
Now using Lemma \ref{Scalc}
\begin{align*}
 \left. \frac{\del}{\del a} S_{j \bk} \right|_{a = 0} =&\ \left.
\frac{\del}{\del a} \left( g(a)^{l \bm} \left( - \del_l \del_{\bm} g(a)_{j \bk}
+ \del g(a) * \delb g(a) \right) \right) \right|_{a = 0}\\
=&\ - h^{l \bm} R_{l \bm j \bk} - g^{l \bm} \del_l \del_{\bm} h_{j \bk}
\end{align*}
Now $- h^{i \bj} R_{i \bj k \bl} = - \overset{\circ}{R}(h)_{k \bl}$ from the
Bianchi identity using that the metric $g(0)$ is K\"ahler-Einstein.
Next, we compute an expression for $\N^* \N h$ using complex coordinates
\begin{align*}
\left(\N^* \N h \right)_{j \bk} =&\ - g^{l \bm} \N_l \N_{\bm} h_{j \bk}\\
=&\ - g^{l \bm} \left( \del_l \del_{\bm} h_{j \bk} - \del_l \gG_{\bm
\bk}^{\bp} h_{j \bp} \right)\\
=&\ - g^{l \bm} \del_l \del_{\bm} h_{j \bk} - R_{\bl}^{\bm} h_{k \bm}\\
=&\ - g^{l \bm} \del_l \del_{\bm} h_{j \bk} - \frac{1}{n} s h_{k \bl}\\
\end{align*}
where $R = S$ is the Ricci tensor of the K\"ahler metric $g(0)$ and $s = r =
\tr_g
S$ is the scalar curvature.  Next we compute using Lemma \ref{sfev}
\begin{align*}
\left. \frac{\del}{\del a} \left( \int_M \tr_g S dV \right) \right|_{a = 0} =&\
\int_M \left( \left< h, - S \right> + \tr_g h \tr_g S \right)
dV\\
=&\ \int_M \left( 1 - \frac{1}{n} \right) \left( \tr_g S \right) \tr_g
h dV\\
=&\ 0
\end{align*}
where the last line follows since $\tr_g S$ is the scalar curvature which is
constant
and $\int_M \tr_g h dV = 0$ since the volume is fixed through $g(s)$.  Thus
\begin{align*}
\left. \frac{\del}{\del a} \frac{1}{n} \left( \int_M \tr_g S dV \right) g
\right|_{a = 0} =&\ \frac{1}{n} \left( \int_M \tr_g S dV \right) h\\
=&\ \frac{1}{n} s h.
\end{align*}
Putting together these calculations yields the result.
\end{proof}
\end{prop}

\begin{defn} Let $L = L(g_0) = \mathcal D_{g_0} \FF$ be the linearization of
$\FF$ at a \static\ metric $g_0$.  We say that $g_0$ is \emph{linearly
stable} if $L \geq 0$.
\end{defn}

\begin{defn} A \static\ metric $g_0$ is \emph{integrable} if for
any solution $h$ of the linearized equation
\begin{gather*}
\mathcal D_{g_0} \left(\FF(g) \right) (h) = 0
\end{gather*}
there exists a path $g(s), s \in (- \ge, \ge)$ of \static\ metrics
where $g(0) = g_0$ and
\begin{gather*}
\frac{d}{ds}_{| s = 0} g(s) = h
\end{gather*}
In particular this implies that the set of Hermitian metrics $g$ satisfying
$\FF(g) \equiv 0$ has a smooth manifold structure near $g_0$.
\end{defn}

We note that by the analysis of Koiso's Theorem it follows that
K\"ahler-Einstein metrics are integrable.  Indeed, any solution to the
linearized deformation equation arises as the variation along a path of
K\"ahler-Einstein metrics, which are static.  This can be seen as follows: If
$c_1(M) < 0$,
K\"ahler-Einstein metrics are linearly stable. If $c_1(M)=0$, any infinitesimal
deformation of
K\"ahler-Einstein metrics is given by Hermitian symmetric deformation of
Einstein metrics which in turn correspond
to $(1,1)$-forms, moreover, we have that
the eigenvalues $L$ are the eigenvalues of the operator
\begin{gather*}
\psi \to \gD_d \psi - \frac{1}{n} s \psi
\end{gather*}
acting on $(1,1)$-forms $\psi$ (\cite{Besse} pg. 362).
If $s = 0$, nonnegativity follows easily and the kernel of $L$ consists of
harmonic $(1,1)$-forms which are
simply variations of K\"ahler-Einstein metrics with vanishing scalar curvature
due to the Calabi-Yau theorem.
So K\"ahler-Einstein metrics are integrable in the case that $c_1(M)=0$.  We now
proceed with the proof of Theorem \ref{Stabthm}.

\begin{proof} Let $(M, g_0, J)$ be a K\"ahler-Einstein manifold.  Fix $h$ a
symmetric two
tensor of type $(1,1)$ such that $\brs{h}_{C^\infty} < \ge' < \ge$ where $\ge'$
and $\ge$ are small positive constants to be chosen later.  We want to show that
solution to the equation
\begin{gather} \label{stabproof10}
\begin{split}
\dt g =&\ - S + Q(T) + \frac{1}{n} \left( \int_M \tr_g \left(S - Q \right) dV
\right) g \\
g(0) =&\ g_0 + h
\end{split}
\end{gather}
exists for all time and converges for $\ge'$ chosen small enough.  Let $h(t) =
g(t) - g_0$.  First consider
\begin{gather} \label{stabproof20}
\begin{split}
\dt h =&\ - S + Q + \frac{1}{n} \left( \int_M \tr_g \left(S - Q \right) dV
\right) g\\
=&\ - \left(\FF(g_0) + \mathcal D \FF_{g_0} (h) + A(g_0, h) \right)\\
=&\ - \mathcal D \FF_{g_0} (h) + A(g_0, h)\\
=&\ - L(h) + A(g_0, h)
\end{split}
\end{gather}
where $A$ represents the higher order terms in the approximation of
$\FF$ by $\mathcal D S_{g_0} (h)$.  Specifically we have
the bounds
\begin{gather} \label{stabproof30}
\brs{A(g_0, h)}_{C^k} \leq C \left( \brs{h}_{C^k} \brs{\N^2 h}_{C^{k-2}} +
\brs{\N h}_{C^{k - 1}}^2 \right)
\end{gather}
where the constant $C$ depends on bounds on the geometry of $g(t)$,
which we are assuming is staying bounded along the flow anyways
since $\brs{g(t) - g_0}_{C^k} < \ge$.  So, fix $T > 0$ and a small
$\ge > 0$.  We would like to show that for $\ge'$ small enough as
above our solution exists on $[0, T)$ and $\brs{h(t)}_{C^k} < \ge$
on this interval.  We start with an $L^2$ growth estimate.

\begin{lemma} \label{stablelemma1} There exists a uniform
(independent of $\ge, \ge', T$) constant $C$ so that if
$\brs{h}_{C^k} < \ge$ for all $t \in [0, T)$, we have
\begin{gather*}
\int_M \brs{h(t)}^2 dV_{g_0} \leq e^{C \ge t} \int_M \brs{h_0}^2
dV_{g_0}
\end{gather*}
\begin{proof} Multiplying the final equation in (\ref{stabproof20})
by $h$ and integrating over $M$ gives
\begin{gather*}
\dt \int_M \brs{h(t)}^2 dV_{g_0} \leq \int_M \left(A * h\right)
dV_{g_0}
\end{gather*}
since $L$ is negative semidefinite.  By straightforward bounds
using integration by parts and the assumed $C^k$ bound on $h$ we are
able to get the bound
\begin{gather*}
\brs{\int_M \left(A * h \right) dV_{g_0}} \leq C \ge \int_M
\brs{h}^2 dV_{g_0}
\end{gather*}
where $C$ depends only on $g_0$.  The result follows immediately.
\end{proof}
\end{lemma}

\begin{lemma} There exists $\ge' = \ge'(T, n) < < \ge$ such that
 if $\brs{h_0}_{C^\infty} < \ge'$, then the solution $g(t)$ exists
 on $[0, T)$ with $\brs{h(t)}_k < \ge$ for all $t \in [0, T)$.
\begin{proof} We use standard parabolic regularity theory.
First we rewrite the evolution equation for $h$ as
\begin{gather} \label{stabproof40}
\dt h = \gD h + \Rm(h) + A(g_0,h)
\end{gather}
Fix a time $\tau < T$.  We first will get an estimate for
$\int_0^\tau \int_M \brs{\N h}^2 dV_{g_0} dt$.  Take the inner
product of (\ref{stabproof40}) with $h$ and integrate over $M$ to get
\begin{gather} \label{stableloc4}
\begin{split}
\frac{1}{2} \dt \int_M \brs{h}^2 =&\ - \int_M \brs{\N h}^2 + \int_M \Rm *
h^{*2} + \int_M \N^2 h * h^{*2} + h * \N h^{*2}\\
\leq&\ - \int_M \brs{\N h}^2 + \gt \int_M \brs{\N h}^2 + C(\gt)
\int_M \brs{h}^2\\
\leq&\ - \frac{1}{2} \int_M \brs{\N h}^2 + C (\gt) \int_M \brs{h}^2
\end{split}
\end{gather}
Using this bound and integrating over time we
conclude
\begin{gather*}
\frac{1}{2} \int_0^\tau \int_M \brs{\N h}^2 \leq \frac{1}{2}
\int_M \brs{h_0}^2 + C(\gt) \tau \sup_{[0, \tau)} \int_M
\brs{h(t)}^2
\end{gather*}
Using Lemma \ref{stablelemma1} we see that $\int_0^\tau \int_M
\brs{\N h}^2$ can be made very small, in particular bounded uniformly in terms
of $\ge'$.  We now show how to get estimates on $\int_0^\tau \int_M \brs{\N^k
h}^2$ for all $k > 0$ in terms of the small constant $\ge'$.
Consider
\begin{align*}
\frac{1}{2} \dt \int_M \brs{\N h}^2 =&\ \int_M \left<\N \left(\gD h + \Rm(h) +
A(g_0, h) \right) , \N h \right> dV\\
&\ + \int_M \N^2 h * \N h * \N h + h * \N h^{*2} dV\\
\leq&\ - \int_M \brs{\N^2 h}^2 + C \int_M \brs{\Rm} \brs{\N h}^2\\
&\ + \theta \int_M \brs{\N^2 h}^2 + C(\theta) \int_M \brs{\N h}^4 + C \ge'
\int_M \brs{\N h}^2 dV\\
\leq&\ - \frac{1}{2} \int_M \brs{\N^2 h}^2 + C \int_M \brs{\N h}^2.
\end{align*}
This implies the bound
\begin{align*}
\frac{1}{2} \int_0^\tau \brs{\N^2 h}^2 \leq&\ \frac{1}{2} \int_M \brs{\N h_0}^2
+ \frac{1}{2} \int_M \brs{h_0}^2 + C \int_0^\tau \int_M  \brs{\N h}^2\\
\leq&\ C \ge'.
\end{align*}
Continuing in this fashion we can induct to get a bound of the above form for
all covariant derivatives of $h$.  Note that for instance we can now bound
\begin{align*}
\int_0^\tau \int_M \brs{\dt h}^2 \leq&\ C \left( \int_0^\tau \int_M \brs{\N^2
h}^2 + \int_0^\tau \int_M \brs{h}^2 \right)\\
\leq&\ C \ge'.
\end{align*}
It is clear that we can in fact get bounds of the form
\begin{gather*}
\int_0^\tau \int_M \brs{\frac{\del^p}{\del t^p} \N^q h}^2 \leq C \ge'
\end{gather*}
for all $p, q > 0$.  One can now apply the Sobolev inequality (with respect to
$g_0$) to conclude $C^k$ bounds on $h$ in terms of $\ge'$.  These bounds will
hold over any time interval where the $L^2$ norm of $h$ is still small.  Since
this time can be made arbitrarily large with small choice of $\ge'$ by Lemma
\ref{stablelemma1} the result follows.
\end{proof}
\end{lemma}

We now improve these estimates to include $L^2$ decay of $h$, which will
ultimately yield the stated long-time existence and convergence.  Say $T$ is a
maximal time such that $\brs{h}_k < \ge$ on $[0, T)$.
Divide the interval $[0, T)$ into intervals of length $\tau$ and let
$N$ be the integer so that $N \tau < T < (N + 1) \tau$.  Let $I_j =
[j \tau, (j + 1) \tau]$.  On $M_j := M \times I_j$ define the inner
product
\begin{gather} \label{stableloc5}
\nm{ f }{M_j} := \int_{j \tau}^{(j + 1) \tau} \nm{f(t)}{L^2(g_0)} dt
\end{gather}
Let $\pi^j$ denote the orthogonal projection onto $\ker \left(\dt +
L \right)$ with respect to $\nm{f}{M_j}$.  Since $L$
is positive semidefinite we see that $\pi$ has no positive
eigenvalues, but there is still the lingering question of zero
eigenvalues.  This is where the integrability property comes
in.  Let $\pi^j_0(h)$ denote the radial component, i.e. the kernel
of $L$. We will show using integrability that there exists a
stationary solution $g_j$ on $M_j$ such that $\pi^j_0 (g(t) - g_j)$ is very
small
compared to $g(t) - g_j$.  This will allow us to conclude $L^2$ decay of $h$ and
then
allow us to conclude convergence.

\begin{lemma} \label{stablelemma2}  Given $\ga > 0$ there exists $\gd = \gd(n,
\tau)$ such that
if
$\sup_{[\tau_0, \tau_0 + \tau]} \brs{h(t)}_k < \gd$ then there exists a
K\"ahler-Einstein
metric $g_1$ such that
\begin{gather} \label{stableloc50}
\brs{\pi^j_0(g - g_1)} \leq \ga (t - \tau_0) \brs{g - g_1}
\end{gather}
and
\begin{gather} \label{stableloc51}
\brs{g_1 - g_0}_{C^k} \leq C \sup_I \brs{g - g_0}_{C^k}
\end{gather}
\begin{proof} Recall from the discussion above that the set of metrics $g$ near
$g_0$ satisfying $\mathcal F(g) = 0$, call it $\mathcal U$, has
a natural smooth manifold
structure.  The tangent space to $\mathcal U$ is given by the kernel
of $L$, call it $\mathcal K$, which is finite dimensional since $L$
is elliptic. Let $\{B_i \}$ be a basis for $\mathcal K$ orthonormal
with respect to the $L^2$ norm induced by $g_0$.  Also using
ellipticity, we get a system of eigenvectors $\{E_\gl \}$ for $L$
orthonormal with respect to the $L^2$ inner product above. We see
that there exist constants $r_\gl$ such that $C_\gl = r_\gl E_\gl
e^{\gl t}$ is a basis for $\ker \left(\dt + L \right)$ which is
orthonormal with respect to the inner product in (\ref{stableloc5}).

Define the map $\Psi : \mathcal U \to \mathcal K$ by $\Psi(g) =
\sum_i \left< g, B_i \right> B_i$.  A simple calculation using the
bases described above shows that for $g_1 \in \mathcal U$,
$\Psi(g_1) = \Psi( \pi^j(g_1)) = \pi^j_0(g_1)$.  Also it is easy to
see that the differential of $\Psi$ at $g_0$ is the identity map, so
we can apply the inverse function theorem.  Fix the time $\tau_0 \in
I_j$.  If $\brs{g(t) - g_0}_k$ is small enough, then in particular
$\pi^j_0(g(\tau_0) - g_0) = \pi^j_0(g(\tau_0))$ can be made small, so that
by the argument above there exists $g_1 \in \mathcal U$ such that
\begin{gather*}
\Psi(g_1) = \pi^j_0(g(\tau_0)).
\end{gather*}
Thus in particular using the above equalities we have $\pi^j_0(g_1 -
g(\tau_0)) = 0$.  Using the evolution equations satisfied by $g$ and
$g_1$ it is clear that one has estimate (\ref{stableloc50}).  Also note that we
have $g_1 = \Psi^{-1}((\pi^j g)_0)$ and
$g_0 = \Psi^{-1}((\pi^j g_0)_0)$ thus using our bound from the
inverse function theorem we get
\begin{gather*}
\nm{g_1 - g_0}{M_j} \leq C \nm{\pi^j(g - g_0)}{M_j}
\end{gather*}
and again using that these are all solutions of the same parabolic
equation, we can get the bound
\begin{gather*}
\brs{g_1 - g_0}_{C^k} \leq C \sup_{I_j} \brs{g - g_0}_{C^k}.
\end{gather*}
\end{proof}
\end{lemma}

\begin{lemma} \label{stablelemma3} Let $I = [\tau_0, \tau_0 + \tau]$
and take $g_1$ as in Lemma \ref{stablelemma2}.  Then there exists $\ge > 0$
depending only on $g_0$ such that if $\brs{h_1(0)}_k < \ge$ where
$h_1 = g - g_1$ then
\begin{gather} \label{stableloc6}
\sup_{\left[ \tau_0 + \frac{\tau}{2}, \tau_0 + \tau \right]} \int_M
\brs{g - g_1}^2 dV_{g_0} \leq e^{- \frac{\tau \gl}{2}} \sup_{\left[
\tau_0, \tau_0 + \frac{\tau}{2} \right]} \int_M \brs{g - g_1}^2
dV_{g_0}
\end{gather}
where $\gl = \min \{ \gl_i : \gl_i \mbox{ is an eigenvalue of } L,
\gl_i \neq 0 \} > 0$.
\begin{proof} Let $h_1(t) = g(t) - g_1$.  If $\brs{h_1(0)}_k < \ge$,
a calculation like that in Lemma \ref{stablelemma1} combined with
the bound on $\pi_0(h_1(t))$ shows that
\begin{align*}
\frac{d}{dt} \int_M \brs{h_1}^2 =&\ \int_M \left< 2 Lh_1, h_1 \right>
dV_{g_0} + \int_M A( h_1, g_0) dV_{g_0}\\
\leq&\ - 2 \gl \int_M \brs{h_1 - \pi^j_0(h_1)} dV_{g_0} + C \ge \int_M
\brs{h_1}^2 dV_{g_0}\\
\leq&\ \left(- \frac{3}{2} \gl + C \ge \right) \int_M \brs{h_1}^2\\
\leq&\ - \gl \int_M \brs{h_1}^2
\end{align*}
as long as $\ge < \frac{\gl}{C}$.  Thus $\int_M \brs{h_1(t)}^2
dV_{g_0} \leq e^{- \gl(t - \tau)} \int_M \brs{h_1(\tau)}^2 dV_{g_0}$
from which the claim follows immediately.
\end{proof}
\end{lemma}

We will need one more Lemma, which roughly says that if a solution
to (\ref{stabproof10}) is decaying at a certain rate at a particular
time then it decayed at that rate earlier in time.  This Lemma is
inspired by Lemma 5.31 in \cite{Tian}, and the proof is the
same.

\begin{lemma} \label{stablelemma4} There exists a constant
$\nu(n, \tau) > 0$ with the following property.  Let $k$ be a
symmetric two-tensor satisfying the equation
\begin{gather*}
\dt k = - L k + A(g_0, k),
\end{gather*}
and
\begin{gather*}
\sup_{ [\tau_0, \tau_0 + \tau] } \brs{k}_{C^k} < \nu
\end{gather*}
and
\begin{gather*}
\brs{\pi_0 (k)} \leq \ga(t - \tau_0) \brs{k}
\end{gather*}
where here we mean projection onto the kernel of $\left( \dt + L
\right)$ restricted to the interval $\left[\tau_0 - \frac{\tau}{2},
\tau_0 + \tau \right]$. Then if
\begin{gather} \sup_{\left[ \tau_0 + \frac{\tau}{2}, \tau_0 + \tau \right]}
\int_M \brs{k}^2 dV_{g_0} \leq e^{- \frac{\tau \gl}{2}} \sup_{\left[
\tau_0, \tau_0 + \frac{\tau}{2} \right]} \int_M \brs{k}^2 dV_{g_0}
\end{gather}
then
\begin{gather}
\sup_{[\tau_0, \tau_0 + \frac{\tau}{2}]} \int_M \brs{k}^2 dV_{g_0}
\leq e^{- \frac{\tau \gl}{2}} \sup_{[\tau_0 - \frac{\tau}{2},
\tau_0]} \int_M \brs{k}^2 dV_{g_0}.
\end{gather}
\begin{proof} First note that the the analogous claim where $k$
satisfies the linear equation $\left( \dt + L \right)k = 0$ and $\pi_0(k) = 0$
is
obvious since $L$ is positive semidefinite and by definition
$\pi(k)_0 = 0$.  In fact there is decay at the rate $\gl$ as opposed
to the $\gl / 2$ in the statement above.  So, if the claim were
false, then for a sequence $\nu_i \to 0$ we would have $k_i$
satisfying the hypothesis but not the conclusion with the bound $\brs{k_i}_{C^k}
<
\nu_i$.  By standard compactness arguments we can parabolically rescale $k$ and
extract a
subsequence converging to $k_\infty$, which satisfies the initial
decay hypothesis but not the conclusion.  Moreover, given that $A$
is quadratic in $k$, it is clear that this $k_\infty$ satisfies the
linear equation $\left( \dt + L \right) k = 0$ and $\pi_0(k) = 0$, contradicting
the
above.
\end{proof}
\end{lemma}

We now proceed with the main proof.  Suppose the maximal existence time
satisfies $T < \infty$, and subdivide $[0, T]$ into $N$
intervals of length $\tau$ labelled $I_j$ as above.  For fixed $j$,
let $g_j$ be the metric such that $\pi^j_0 (g(t) - g_j) = 0$ on
$I_j$ given by Lemma \ref{stablelemma2}.  Define $h_j := g(t) -
g_j$.  By Lemma \ref{stablelemma3} and parabolic regularity we have
that
\begin{gather}
\sup_{\left[ \left(j + \frac{1}{2} \right) \tau, \left(j + 1 \right)
\tau \right]} \brs{h_j} \leq C e^{- \frac{\tau \gl}{2}} \sup_{\left[
j \tau, \left(j + \frac{1}{2} \right) \tau \right]} \brs{h_j}.
\end{gather}
And we can apply Lemma \ref{stablelemma4} inductively to conclude
\begin{gather*}
\sup_{[j \tau, (j + 1) \tau]} \brs{h_j} \leq e^{- \gl \tau(j - 1)}
\sup_{\left[0, \frac{\tau}{2} \right]} \brs{h_j}.
\end{gather*}
This allows us to conclude that on $I_j$ we have
\begin{align*}
\brs{\dt g} =&\ \brs{\dt \left(g - g_j \right)}\\
\leq&\ C \sup_{I_j} \brs{h_j}_k\\
\leq&\ C \ge e^{- \gl \tau (j - 1)} = \frac{C \ge}{p^{j - 1}}.
\end{align*}
Now note that simply integrating over time we see that
\begin{gather*}
\sup_{I_j} \brs{g - g_0} \leq 2 \tau \sup_{I_j \cup I_{j - 1}}
\brs{\dt g} + \sup_{I_{j-1}} \brs{g - g_0}.
\end{gather*}
Applying this estimate inductively we see that
\begin{gather} \label{stableloc7}
\begin{split}
\sup_{I_j} \brs{g - g_0} \leq&\ 2 \tau \sum_{k=1}^{N} \sup_{I_k \cup
\dots \cup I_N} \brs{\dt g} + \sup_{I_0} \brs{g - g_0}\\
\leq&\ \sum_{k = 1}^\infty \frac{2 \tau C \ge}{p^{k - 1}} + \sup_{I_0}
\brs{h_0}\\
\leq&\ \frac{2 \tau C \ge}{p - 1} + \sup_{I_0} \brs{h_0}.
\end{split}
\end{gather}
Now we want to choose our constants $\ge, \ge'$ and $\tau$ to derive
a contradiction from this equation.  So, choose $\tau$ initially so
large that
\begin{gather}
\frac{1}{c(n) e^{\tau \gl}} + \frac{2 C \tau}{e^{\tau \gl} - 1} <
\frac{1}{C} e^{- \frac{\tau \gl}{4}}
\end{gather}
where $c(n)$ is a fixed large constant and $C$ is a constant
depending only on $g_0$.

Now let $\ge = \min \{ \gd(n, \tau),\nu(n,\tau),\frac{\gl}{C_0} \}$
where $\gd(n, \tau)$ is as in Lemma \ref{stablelemma2},
$\nu(n,\tau)$ is as in Lemma \ref{stablelemma4}, and $C_0$ is a
constant depending only on $g_0$ and the dimension which we now make
explicit.  By Lemma \ref{stablelemma1} we can bound the growth of
$L^2$ derivatives of $h$, and then by Sobolev embeddings we can
bound $C^k$ norms. Specifically there exists a constant depending
only on $g_0$ so that
\begin{gather*}
\brs{h}_{C^4} < C e^{C \ge \tau} \brs{h(0)}_{C^0}
\end{gather*}
Then let $C_0 := 12 C$.  Note that if $\ge < \frac{\gl}{C_0}$ and we
start our flow with some $h(t_0)$ satisfying $\brs{h(t_0)}_{C^4} <
\frac{\ge}{C} e^{- \tau \gl/4}$, the solution exists at least on
$[t_0, t_0 + 3 \tau)$ and moreover $\sup_{[t_0, t_0 + 3 \tau)}
\brs{h(t)}_{C^4} < \ge$.

Now again using Lemma \ref{stablelemma1} we see that we may choose
$\ge'$ so that the solution exists on $[0, 3 \tau]$ and further
\begin{gather*}
\sup_{[0, 3 \tau]} \brs{h(t)}_k < \frac{\ge}{c(n)} e^{- \tau \gl}.
\end{gather*}
Since $\ge < \min \{ \gd(n, \tau), \nu(n, \tau) \}$ we can apply
(\ref{stableloc7}) to get that
\begin{align*}
\sup_{I_N} \brs{g - g_0} \leq&\ \ge \left( \frac{1}{c(n)
e^{\tau \gl}} + \frac{2 C \tau}{e^{\tau \gl} - 1} \right)\\
\leq&\ \frac{\ge}{C} e^{- \frac{\tau \gl}{4}}.
\end{align*}
Thus the solution to HCF with initial metric $g(T - \tau)$ exists on an interval
of length $3 \tau$ with $\brs{h}_{C^2} < \ge$, contradicting the
maximality of $T$. Thus the solution exists for all time and
$\brs{g(t) - g_0}_k < \ge$ for all time.  Indeed we have decay
\begin{gather*}
\brs{g(t) - g_j} \leq C e^{- \gl t}
\end{gather*}
for all $t \in [0, j \tau)$ and for all $j$.  Since $\{ g_j \}$ is a
sequence of K\"ahler-Einstein metrics with uniform $C^k$ bounds, we get a
convergent subsequence $g_j \to g_\infty$ a critical metric, with
exponential convergence $g(t) \to g_\infty$.
\end{proof}

\section{Further Questions}

It bears mentioning that the Theorem 1.1 and Theorem 1.2 are both true for more
general equations.  Specifically both results hold for solutions to (\ref{HCF})
where $Q$ can be \emph{any} quadratic expression in the torsion.

The HCF is similar in some regards to certain renormalization group flows
arising in physics where external fields, say Yang-Mills or B-fields, are added
to the pure gravity theory and then arise in the flow equations, see for
instance \cite{OSW}, \cite{Streets1}, \cite{Streets2}.  In these flows the
torsion is given as an
external field, whereas in HCF everything is defined in terms of the metric.  A
similar case is studied in \cite{Bryant}, \cite{Hitchin} where a ``holonomy
flow'' is proposed for closed
$G_2$ structures.  Here one evolves the definite three-form $\sigma$ defining
the $G_2$ structure by the Hodge Laplacian of $\sigma$ taken with the metric
induced by $\sigma$.  This is a quasilinear equation which bears a certain
resemblance to HCF in that it can be written as ``Ricci flow plus torsion,''
where the torsion is defined in terms of the underlying metric.  The techniques
of our stability theorem likely apply to show stability of this flow near
$G_2$-holonomy spaces with negative semidefinite Lichnerowicz operator.

Finally, Hermitian curvature flow provides a framework for addressing questions
on the existence of
integrable complex structures.  In particular, if one had a complete description
of the behavior of this flow for certain geometric conditions and a complete
understanding of the limiting objects, one could then describe the manifolds
admitting integrable complex structures with Hermitian metrics satisfying the
initial geometric conditions.  With strong enough convergence results for this
flow one could in particular answer the question of the existence of an
integrable complex structure on $S^6$.  Since we expect our flow to be moving
towards K\"ahler metrics and we know that $S^6$ does not support K\"ahler
metrics, we expect the volume
normalized flow to develop singularities either at finite time or at infinity.
In either case these singularities will have some extra structure and can
possibly be classified. It may lead to a contradiction to the assumption on the
existence of integrable complex structures on $S^6$.

\section{Appendix: Variational Formulas}
In this appendix we collect variational formulas for quantities related to the
curvature and torsion of Hermitian metrics.

\begin{lemma} \label{connectionvariation} Let $g(a)$ be a family of Hermitian
metrics compatible with the
given complex structure $J$.  Then
\begin{align*}
\frac{\del}{\del a} \gG_{i k}^l =&\ g^{l \bm} \N_i h_{k \bm}\\
\frac{\del}{\del a} \gG_{\bi \bk}^{\bl} =&\ g^{\bl m} \N_{\bi} h_{m \bk}
\end{align*}
\begin{proof} We compute directly in canonical complex coordinates
\begin{align*}
\frac{\del}{\del a} \gG_{i k}^l =&\ \frac{\del}{\del a} g^{l \bm} \left( \del_i
g_{k \bm} \right)\\
=&\ - h^{l \bm} \del_i g_{k \bm} + g^{l \bm} \del_i h_{k \bm}\\
=&\ g^{l \bm} \left( \del_i h_{k \bm} - \gG_{i k}^p h_{p \bm} \right)\\
=&\ g^{l \bm} \N_i h_{k \bm}.
\end{align*}
Which gives the first formula and the second follows by conjugation.
\end{proof}
\end{lemma}

\begin{lemma} \label{curvaturevariation} Let $g(a)$ be a family of Hermitian
metrics compatible with the
given complex structure $J$.  Then
\begin{align*}
\frac{\del}{\del a} \Omega_{i \bj k}^l =&\ - g^{m \bl} \N_{\bj} \N_i h_{k \bm}\\
\frac{\del}{\del a} \Omega_{i \bj k \bl} =&\ \Omega_{i \bj k}^m h_{m \bl} -
\N_{\bj} \N_i h_{k \bl}
\end{align*}
\begin{proof} We compute directly
\begin{align*}
\frac{\del}{\del a} \Omega_{i \bj k}^l =&\ \frac{\del}{\del a} \left( -
\del_{\bj} \gG_{i k}^l \right)\\
=&\ - \del_{\bj} \left( g^{l \bm} \N_i h_{k \bm} \right)\\
=&\  g^{l \bp} \del_{\bj} g_{\bp q} g^{q \bm} \N_i h_{k \bm} - g^{l \bm}
\del_{\bj} \N_i h_{k \bm}\\
=&\ - g^{m \bl} \left( \del_{\bj} \N_i h_{k \bm} - \gG_{\bj \bm}^{\bp} \N_i h_{k
\bp} \right)\\
=&\ - g^{m \bl} \N_{\bj} \N_i h_{k \bm}.
\end{align*}
This gives the first formula and the second follows easily.
\end{proof}
\end{lemma}

\begin{lemma} \label{svar} Let $g(a)$ be a family of Hermitian metrics
compatible with the
given complex structure $J$.  Then
\begin{align*}
\frac{\del}{\del a} s =&\ - \gD \tr h - \left< h, S + \divg^{\N} T - \N w
\right>.
\end{align*}
\begin{proof}
We compute using Lemma \ref{curvaturevariation}
\begin{align*}
\frac{\del}{\del a} s =&\ \frac{\del}{\del a} g^{k \bl} S_{k \bl}\\
=&\ - h^{k \bl} S_{k \bl} + g^{k \bl} \left[ - \gD h_{k \bl} - h^{i \bj}
\Omega_{i \bj k \bl} + h_{k \bm} S^{\bm}_{\bl} \right]\\
=&\ - \gD \tr h - \left< h, P \right>.
\end{align*}
The result now follows from Lemma \ref{PSForm}.
\end{proof}
\end{lemma}

\begin{lemma} \label{torsionvariation} Let $g(a)$ be a family of Hermitian
metrics compatible with the
given complex structure $J$.  Then
\begin{align*}
\frac{\del}{\del a} T_{i j \bk} =&\ \N_i h_{j \bk} - \N_j h_{i \bk} + T_{i j}^m
h_{m \bk}\\
\frac{\del}{\del a} w =&\ \N \tr h - \divg^{\N} h.
\end{align*}
\begin{proof} We compute directly
\begin{align*}
\frac{\del}{\del a} T_{i j \bk} =&\ \frac{\del}{\del a} \left( \del_i g_{j \bk}
- \del_j g_{i \bk} \right)\\
=&\ \del_i h_{j \bk} - \del_j h_{i \bk}\\
=&\ \N_i h_{j \bk} + \gG_{i j}^m h_{m \bk} - \N_j h_{i \bk} - \gG_{j i}^p h_{p
\bk}\\
=&\ \N_i h_{j \bk} - \N_j h_{i \bk} + T_{i j}^m h_{m \bk}.
\end{align*}
This gives the first formula.  For the second we see
\begin{align*}
\frac{\del}{\del a} w_i =&\ \frac{\del}{\del a} g^{j \bk} T_{i j \bk}\\
=&\ - h^{j \bk} T_{i j \bk} + g^{j \bk} \left(\N_i h_{j \bk} - \N_j h_{i \bk} +
T_{i j}^m h_{m \bk} \right)\\
=&\ \N_i \tr h - \divg^{\N} h_i.
\end{align*}
\end{proof}
\end{lemma}

\begin{lemma} \label{q1variation} Let $g(s)$ be a family of Hermitian metrics
compatible with the
given complex structure $J$.  Then
\begin{align*}
\frac{\del}{\del s} \brs{T}^2 =&\ \left< h, -2 Q^1 + Q^2 \right> + 4 \left<\N h,
T \right>
\end{align*}
where
\begin{align*}
\left<\N h, T \right> = \frac{1}{2} g^{i \bj} g^{k \bl} g^{m \bn} \left(\N_i
h_{k \bn} T_{\bj \bl m} + T_{i k \bn} \N_{\bj} h_{\bl m} \right).
\end{align*}
\begin{proof} We compute directly
\begin{align*}
\frac{\del}{\del s} \brs{T}^2 =&\ \frac{\del}{\del s} g^{i \bp} g^{j \bq} g^{\bk
r} T_{i j \bk} T_{\bp \bq r}\\
=&\ - h^{i \bp} g^{j \bq} g^{\bk r} T_{i j \bk} T_{\bp \bq r} - h^{j \bq} g^{i
\bp} g^{\bk r} T_{i j \bk} T_{\bp \bq r} - h^{\bk r} g^{i \bp} g^{j \bq} T_{i j
\bk} T_{\bp \bq r}\\
&\ + g^{i \bp} g^{j \bq} g^{\bk r} \left[\left( \N_i h_{j \bk} - \N_j h_{i \bk}
+ T_{i j}^m h_{m \bk} \right) T_{\bp \bq r} + T_{i j \bk} \left(\N_{\bp} h_{\bq
r} - \N_{\bq} h_{\bp r} + T_{\bp \bq}^{\bs} h_{r \bs} \right) \right]\\
=&\ \left<h, -2 Q^1 + Q^2 \right> + 4 \left<\N h, T \right>.
\end{align*}
The result follows.
\end{proof}
\end{lemma}

\begin{lemma} \label{q3variation}
Let $g(s)$ be a family of Hermitian metrics compatible with the given complex
structure $J$.  Then
\begin{align*}
\frac{\del}{\del s} \brs{w}^2 =&\ - \left< h, Q^3 \right> + 2 \left< \N \tr h -
\divg^{\N}
h, w \right>.
\end{align*}
\begin{proof} In canonical complex coordinates at a point we compute
\begin{align*}
\frac{\del}{\del s} &g^{i \bj} g^{m \bn} g^{r \bs} T_{i m \bn} T_{\bj \bs r}\\
=&\ - h^{i \bj} g^{m \bn} g^{r \bs} T_{i m \bn} T_{\bj \bs r} - g^{i \bj} h^{m
\bn} g^{r \bs} T_{i m \bn} T_{\bj \bs r} - g^{i \bj} g^{m \bn} h^{r \bs} T_{i m
\bn} T_{\bj \bs r}\\
&\qquad + g^{i \bj} g^{m \bn} g^{r \bs} \left[ \left( \N_i h_{m \bn} - \N_m
h_{i \bn} + T_{i m}^{p} h_{p \bn} \right) T_{\bj \bs r} \right.\\
& \qquad \left. + \left( \N_{\bj} h_{\bs r} - \N_{\bs}
h_{\bj r} + T_{\bj \bs}^{\bq} h_{\bq r} \right) T_{i m \bn} \right]\\
=&\ - \left< h, Q^3 \right> + 2 \left< \N \tr h - \divg^{\N} h, w \right>.
\end{align*}
\end{proof}
\end{lemma}

\begin{lemma} \label{sfev} Let $g(a)$ be a family of Hermitian metrics
compatible with the given complex structure $J$.  Then
\begin{align*}
\frac{\del}{\del a} \int_M s dV =&\ \int_M \left[ \left< h, - S - \divg^{\N} T +
\N w \right> + \tr h \left( s - \divg^{\N} w - \brs{w}^2 \right) \right] dV.
\end{align*}
\begin{proof} From Lemma \ref{svar} we see
\begin{align*}
\frac{\del}{\del a} \int_M s dV =&\ \int_M \left[ - \gD \tr h + \left< h, - S -
\divg^{\N} T + \N w \right> + s \tr h \right] dV.
\end{align*}
Now by Lemma \ref{ibplemma3} we see
\begin{align*}
\int_M - \gD \tr h dV =&\ \int_M \tr h \left( - \divg^{\N} w - \brs{w}^2 \right)
dV
\end{align*}
and the result follows.
\end{proof}
\end{lemma}

\begin{lemma} \label{q1fev} Let $g(a)$ be a family of Hermitian metrics
compatible with the
given complex structure $J$.  Then
\begin{align*}
\frac{\del}{\del a} \int_M \brs{T}^2 dV =&\ \int_M \left[ \left< h, - 2 Q^1 +
Q^2 - 4 Q^4 - 4 \divg^{\N} T \right> \right.\\
&\ \qquad \quad + \left. \left(\tr h \right) \brs{T}^2 \right] dV.
\end{align*}
\begin{proof} We directly compute using Lemma \ref{q1variation}
\begin{align*}
\frac{\del}{\del a} \int_M \brs{T}^2 dV =&\ \int_M \left[ \left< h, - 2 Q^1 +
Q^2  \right> + 4 \left< \N  h , T \right> + \left(\tr h \right) \brs{T}^2
\right] dV.
\end{align*}  We now integrate by parts
\begin{align*}
\int_M \left< \N h, T \right> dV =&\ \int_M g^{i \bp} g^{j \bq} g^{\bk r}
\N_i h_{j \bk} T_{\bp \bq r} dV\\
=&\ \int_M g^{i \bp} g^{j \bq} g^{\bk r} \left(\del_i h_{j \bk} - \gG_{i j}^{s}
h_{s \bk} \right) T_{\bp \bq r} dV\\
=&\ - \int_M h_{j \bk} \del_i \left( g^{i \bp} g^{j \bq} g^{\bk r} T_{\bp \bq r}
dV \right) + g^{i \bp} g^{j \bq} g^{\bk r} \gG_{i j}^{s} h_{s \bk} T_{\bp
\bq r} dV\\
=&\ \int_M h_{j \bk} \left[ g^{i \bn} \del_i g_{m \bn} g^{m \bp} g^{j \bq}
g^{\bk r} T_{\bp \bq r} + g^{i \bp} g^{j \bn} \del_i g_{m \bn} g^{m \bq} g^{\bk
r} T_{\bp \bq r} \right.\\
&\ \left. \qquad \quad + g^{i \bp} g^{j \bq} g^{\bk m} \del_i g_{m \bn} g^{\bn
r} T_{\bp \bq r} - g^{i \bp} g^{j \bq} g^{\bk r} \del_i T_{\bp \bq r} \right.\\
&\ \left. \qquad \quad - g^{i \bp} g^{j \bq} g^{\bk r} T_{\bp \bq r} g^{m \bn}
\del_i g_{m \bn} \right.\\
&\ \left. \qquad \quad - g^{m \bp} g^{n \bq} g^{\bk r} \gG_{m n}^j T_{\bp \bq r}
\right] dV\\
=&\ \int_M \left< h, - \divg^{\N} T - Q^4 \right> dV
\end{align*}
Plugging in this simplification yields the result.
\end{proof}
\end{lemma}

\begin{lemma} \label{q3fev} Let $g(a)$ be a family of Hermitian metrics
compatible with the
given complex structure $J$.  Then
\begin{align*}
\frac{\del}{\del a} \int_M \brs{w}^2 dV =&\ \int_M \left[ \left< h, Q^3 + 2
\N w \right> + \tr h \left[ - 2 \divg^{\N} w - \brs{w}^2 \right] \right] dV.
\end{align*}
\begin{proof} We directly compute using Lemma \ref{q3variation} and integrating
by
parts using Lemmas \ref{ibplemma2} and \ref{ibplemma4}
\begin{align*}
\frac{\del}{\del a} \int_M \brs{w}^2 dV =&\ \int_M \left[ \left< h, - Q^3
\right> + 2 \left< \N \tr h - \divg^{\N} h, w \right> + \left(\tr h \right)
\brs{w}^2 \right] dV\\
=&\ \int_M \left[ \left< h, - Q^3 \right> + 2 \tr h \left[ - \divg^{\N} w -
\brs{w}^2 \right] \right.\\
&\ \left. \qquad \quad + 2 \left< h, \N w + Q^3
\right> + (\tr h) \brs{w}^2\right] dV
\end{align*}
\end{proof}
\end{lemma}

\begin{lemma} \label{ibplemma2} Given $\phi \in C^{\infty}(M)$ and $\ga \in
T^{1,0}(M)$ we have
\begin{align*}
\int_M \left< \N \phi, \ga \right> =&\ \int_M \phi \left[ - \divg^{\N} \ga -
\left< w, \ga \right> \right]
dV.
\end{align*}
\begin{proof} We directly compute
\begin{align*}
\int_M \left< \N \phi, \ga \right> dV =&\ \int_M g^{i \bj} \del_i \phi
\ga_{\bj} dV\\
=&\ - \int_M \phi \del_i \left[ g^{i \bj} \ga_{\bj} dV \right]\\
=&\ \int_M \phi \left[ g^{i \bs} \del_i g_{r \bs} g^{r \bj} \ga_{\bj} - g^{i
\bj} \del_i \ga_{\bj} - g^{i \bj} \ga_{\bj} g^{p
\bq} \del_i g_{p \bq} \right] dV\\
=&\ \int_M \phi \left[ - \divg^{\N} \ga - \left< w, \ga \right> \right]
dV.
\end{align*}
\end{proof}
\end{lemma}

\begin{lemma} \label{ibplemma3} Given $\phi \in C^{\infty}(M)$ we have
\begin{align*}
\int_M \gD \phi dV =&\ \int_M \phi \left[ \divg^{\N} w + \brs{w}^2 \right] dV
\end{align*}
\begin{proof}
We directly compute
\begin{align*}
\int_M \gD \phi dV =&\ \int_M g^{i \bj} \del_{\bj} \del_i \phi dV\\
=&\ - \int_M \del_i \phi \del_{\bj} \left( g^{i \bj} dV \right)\\
=&\ \int_M \del_i \phi \left[ g^{i \bs} \del_{\bj} g_{r
\bs} g^{r \bj} - g^{i \bj} g^{p \bq} \del_{\bj} g_{p \bq} \right] dV\\
=&\ \int_M \del_i \phi \left[ g^{i \bj} g^{p \bq} T_{\bq \bj p} \right] dV\\
=&\ - \int_M \left< \N \phi, w \right> dV\\
=&\ \int_M \phi \left[ \divg^{\N} w + \brs{w}^2 \right] dV.
\end{align*}
In the line we applied Lemma \ref{ibplemma2}.
\end{proof}
\end{lemma}

\begin{lemma} \label{ibplemma4} Given $\gb \in T^{1,0} M$
and $h \in \Sym^{1,1}(M)$ we have
\begin{align*}
\int_M \left< \divg^{\N} h, \gb \right> dV =&\ \int_M \left[ \left< h, -
\N \gb - w \otimes \gb \right>
\right] dV
\end{align*}
\begin{proof}
\begin{align*}
\int_M \left< \divg^{\N} h, \gb \right> dV =&\ \int_M g^{\bj l} g^{\bk i}
\N_{\bk} h_{i \bj} \gb_l dV\\
=&\ \int_M g^{i \bk} g^{l \bj} \left[ \del_{\bk} h_{i \bj} - \gG_{\bk \bj}^{\bp}
h_{i \bp} \right] \gb_l dV\\
=&\ - \int_M h_{i \bj} \del_{\bk} \left[ g^{i \bk} g^{l \bj} \gb_l dV
\right] + g^{i \bk} g^{l \bj} \gG_{\bk \bj}^{\bp} h_{i \bp} \gb_l dV\\
=&\ \int_M h_{i \bj} \left[ g^{i \bs} \del_{\bk} g_{r \bs} g^{r \bk} g^{l \bj}
\gb_l + g^{i \bk} g^{l \bs} \del_{\bk} g_{r \bs} g^{r \bj} \gb_l - g^{i \bk}
g^{l \bj} \del_{\bk} \gb_l \right.\\
&\ \left. \quad \qquad - g^{i \bk} g^{l \bj} g^{r \bs} \del_{\bk} g_{r \bs}
\gb_l \right] dV - g^{i \bk} g^{l \bj} \gG_{\bk \bj}^{\bp} h_{i \bp}
\gb_l dV\\
=&\ \int_M \left[ \left< h, - \N \gb - w \otimes \gb \right>
\right] dV.
\end{align*}
\end{proof}
\end{lemma}

\bibliographystyle{hamsplain}

\end{document}